\documentclass[final,leqno]{siamltex}

\usepackage{amsmath,amssymb}
\usepackage[ruled,vlined,noend]{algorithm2e}
\usepackage{graphicx}\usepackage{dcolumn}\usepackage{booktabs}
\usepackage{bm}\usepackage[font=footnotesize,format=plain,labelfont=bf]{caption}
\usepackage{subcaption}
\usepackage{pstricks}
\usepackage{psfrag,rotating,multirow,setspace,colordvi}
\usepackage{xspace}

\newcommand{\bbf}{\mathbf{f}}
\newcommand{\bc}{\mathbf{c}}

\newcommand{\bn}{\mathbf{n}}

\newcommand{\bx}{\mathbf{x}}

\newcommand{\bR}{\mathbf{R}}
\newcommand{\br}{\mathbf{r}}

\newcommand{\bxi}{\boldsymbol{\xi}}
\newcommand{\bd}{\boldsymbol{\delta}}

\newcommand{\mC}{\mathcal{C}}

\newcommand{\mO}{\mathcal{O}}
\newcommand{\mQ}{\mathcal{Q}}

\newcommand{\mW}{\mathcal{W}}

\newcommand{\dx}{\mathrm{d}\bx}

\newcommand{\erfc}{\mathrm{erfc}}

\newcommand{\phsr}{\Phi^{\text{sr}}}
\newcommand{\phsm}{\Phi^{\text{sm}}}
\newcommand{\phsc}{\Phi^{\text{sc}}}

\newcommand{\nel}{N_\text{el}}
\newcommand{\nelsr}{N_\text{el}^\text{sr}}
\newcommand{\nelx}{n_\text{el}^x}
\newcommand{\nely}{n_\text{el}^y}
\newcommand{\nelz}{n_\text{el}^z}
\newcommand{\nsctotal}{N_\text{sc}}
\newcommand{\nsc}{n_\text{sc}}

\newcommand{\scbasis}{\psi}
\newcommand{\bscbasis}{\bar{\psi}}
\newcommand{\mbasis}{u^\text{m}}
\newcommand{\nmesh}{N_\text{m}}
\newcommand{\rhom}{\rho_\text{m}}
\newcommand{\bxmesh}{\bx^{\text{m}}}
\newcommand{\bxcharge}{\bx^{\text{c}}}

\newcommand{\pppm}{P$^3$M\xspace}

\newtheorem{remark}{Remark}

\usepackage[colorinlistoftodos, color=orange, bordercolor=orange]{todonotes}

\title{A Finite Element Based \pppm\ Method for $N$-body Problems}

\author{Natalie N.~Beams\footnotemark[1]
\and
Luke N.~Olson\footnotemark[2]
\and
Jonathan B.~Freund\footnotemark[3]
}

\begin{document}

\maketitle
\renewcommand{\thefootnote}{\fnsymbol{footnote}}
\footnotetext[1]{Department of Mechanical Science \& Engineering,
                 University of Illinois at Urbana-Champaign, Urbana, IL 61801,
                 \texttt{beams2@illinois.edu}}
\footnotetext[2]{Department of Computer Science,
                 University of Illinois at Urbana-Champaign, Urbana, IL 61801,
                 \texttt{lukeo@illinois.edu}}
\footnotetext[3]{Department of Mechanical Science \& Engineering, and
                 Department of Aerospace Engineering,
                 University of Illinois at Urbana-Champaign, Urbana, IL 61801,
                 \texttt{jbfreund@illinois.edu}}
\renewcommand{\thefootnote}{\arabic{footnote}}

\begin{abstract}
We introduce a fast mesh-based method for computing $N$-body interactions that
is both scalable and accurate.  The method is founded on a
particle-particle{--}particle-mesh (\pppm) approach, which decomposes a potential
into rapidly decaying short-range interactions and smooth, mesh-resolvable
long-range interactions.  However, in contrast to the traditional approach of
using Gaussian screen functions to accomplish this decomposition, our method
employs specially designed polynomial bases to construct the screened
potentials.  Because of this form of the screen, the long-range component of the
potential is then solved \textit{exactly} with a finite element method, leading
ultimately to a sparse matrix problem that is solved efficiently with standard
multigrid methods.  Moreover, since this system represents an exact
discretization, the optimal resolution properties of the FFT are unnecessary,
though the short-range calculation is now more involved than \pppm/PME methods.
We introduce the method, analyze its key properties, and demonstrate the
accuracy of the algorithm.
\end{abstract}

\begin{keywords}
$N$-body, finite element, multigrid, \pppm, PME, multipole methods
\end{keywords}

\begin{AMS}
70--08, 70F10, 65N30, 65N99
\end{AMS}

\pagestyle{myheadings}
\thispagestyle{plain}
\markboth{N.~N.~Beams, L.~N.~Olson, and J.~B.~Freund}{FE-based $N$-body Method}

\section{Introduction}

$N$-body interactions arise in a range
of applications, including molecular dynamics, plasma dynamics, vortex
methods, and viscous flow: systems that are described by
a Green's function solution to the Poisson equation or its
derivatives.  We focus on three-dimensional electrostatic-like $1/R$
interactions, where $R$ is the distance to a particle; this is the
simplest kernel in three dimensions and well-known for
this class of problems.  However, the resulting algorithm we describe
extends to other systems.  We consider a periodic domain,
which is commonly used to model extensive systems, and discuss a
straightforward extension to other boundary conditions in Section~\ref{sec:bc}.
Without loss of generality we consider a $L^3$ cubic unit cell containing $N$
point charges, which has the total electrostatic potential energy
\begin{align}
  \mathcal{U} & = \frac{1}{2}\sum_{\bn =
-\boldsymbol{\infty}}^{\boldsymbol{\infty}}
\sum_{\substack{i,j=1\\ i\ne j, \bn = 0}}^{N}
\frac{Q_i Q_j}{|\bx_i - \bx_j + \bn L|} \equiv \frac{1}{2} \sum_{i=1}^{N} Q_i \Phi_i,\label{eq:potentialenergy}
\end{align}
where $\Phi_i$ is the electrostatic potential at location $\bx_i$ of
particle $i$ with charge $Q_i$.   The central challenge in
(\ref{eq:potentialenergy}) is the computation of the
potential,
\begin{equation}
 \Phi(\bx_i) = \sum_{\bn =
-\boldsymbol{\infty}}^{\boldsymbol{\infty}}
\sum_{\substack{j=1\\ i\ne j, \bn = 0}}^{N}
\frac{Q_j}{|\bx_i - \bx_j + \bn L|},\label{eq:potential}
\end{equation}
because of the fairly slow $1/R$ decay rate of the interactions at large
distances.

There are a number of approaches for efficiently evaluating
(\ref{eq:potential}).  The most widely used methods are generally classified as
either tree-based, such as the fast multipole method
(FMM)~\cite{Greengard:1987}, or mesh-based (sometimes called
``particle-in-cell''), such as the particle-particle{--}particle-mesh (\pppm)
method~\cite{Hockney:1988} and its popular variant, the particle-mesh-Ewald
(PME) method~\cite{Darden:1993, Essmann:1995}.  In the FMM, particles are
grouped within multipole expansions to provide an accurate representation of
their combined influence at a distance, thus limiting the number of terms
needed to explicitly compute the interactions.  The resulting algorithm
scales with $\mO(N)$ complexity, although the coefficient in this scaling can be
large, especially if a high-order multipole expansion is required for the
desired accuracy~\cite{Greengard:1997}.  Efficient implementations are
intricate---especially in parallel---but demonstrated, and the FMM has been
shown to be effective as an adaptive three-dimensional
algorithm~\cite{Cheng:1999}.  The method also extends to systems with more
complicated kernels, such as Stokes flow~\cite{Tornberg:2008,Wang:2007,Veerapaneni:2011}.

In comparison, mesh-based methods also reduce the number of explicit
calculations but achieve this by splitting the potential into a rapidly
decaying component $\phsr$, which is accurately calculated with inclusion
of only a few short-range interactions, and a smooth part $\phsm$, which is
solved on a mesh covering the domain~\cite{Hockney:1988}.  It is
instructive to view this splitting as the addition and subtraction of
strategically selected ``screening'' functions, so that the potential in
(\ref{eq:potential}) decomposes as
\begin{equation}
\label{eq:phisplit}
\Phi_i = \underbrace{\Phi_i - \phsm_i}_{\phsr_i} + \phsm_i.
\end{equation}

\begin{figure}[!ht]
 \centering
 \includegraphics[width=0.8\textwidth]{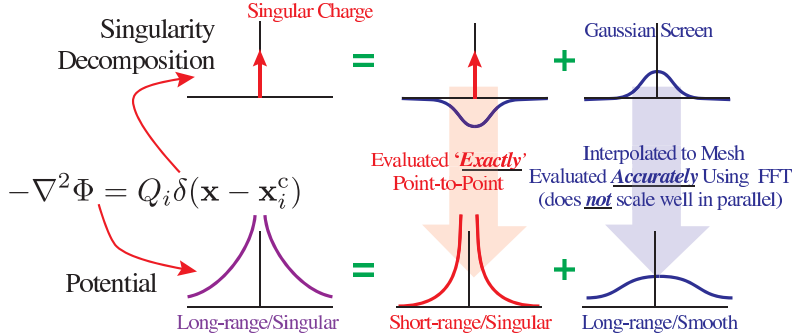}
   \caption{Introduction of a Gaussian screen to define (\ref{eq:phisplit}) and resulting decomposition of the potential of the singular charge.}\label{fig:splitting}
\end{figure}

The particle-mesh-Ewald (PME) method~\cite{Darden:1993} bases this decomposition
directly on the Ewald summation~\cite{Ewald:1921} for (\ref{eq:potential}) and
uses Lagrangian interpolation to move between particle locations and the mesh,
while the smooth PME (SPME) uses B-spline interpolants, similar to those
proposed in the P$^3$M method~\cite{Essmann:1995}.  PME-based algorithms use
Gaussian screening functions, as illustrated in Figure~\ref{fig:splitting}.
Here, the screen is designed to yield a $\phsr$ that is straightforward to
calculate within a prescribed cutoff at radius $R_c$, while the long-range
portion of the potential remains smooth.
In PME-based methods, the Gaussian screen yields a $\phsm$ that is accurately
solved by fast Fourier transforms (FFTs).  To do this, the screen is
interpolated to a regular mesh and the Poisson (or similar) operator is
inverted.  For computational efficiency it is desirable that these screens be
as compact as possible and barely resolved on the mesh, since this maximizes
the decay of the screened potential $\phsr$.  Fast decay allows for a small
point-to-point interaction cutoff distance $R_c$, which reduces the number of
interactions that need to be explicitly computed for the targeted accuracy.
The ideal wavenumber resolution of the FFT provides accurate representation of
the most compact screens possible.    The FFT also makes these methods most
natural for periodic domains, but they can be extended to free space
\cite{Hockney:1988, Pollock:1996, Freund:2002}.

Here we propose a fundamentally new decomposition that is constructed within a
\pppm-type framework.  The method incorporates screen potentials that are
selected to yield exact mesh potentials, which has many potential benefits.  The
screens are designed for a mesh and thus have no explicit dependence on problem
geometry; this suggests complex geometries as well as more general boundary
conditions fit naturally within this method.  In addition, the exact mesh
potential recast the problem as sparse matrix problem where multigrid methods
are known (and shown in Section 4) to be effective and scale to high core
counts~\cite{2012_BaAl_etal_Scalability}.  As a result, since the method does
not rely on the Fourier resolution for an accurate mesh solution, a global FFT
can be avoided, which may be beneficial at extreme scales.  Indeed, while
multigrid methods are ultimately latency bound, they do not exhibit the strong
dependence on a machine's half-bandwidth, which is a limited factor of using
multidimensional FFTs a large core counts~\cite{2010_GhGr_feasibility}.

The new decomposition we propose comes with the cost of representing more
intricate short-range interactions.  The calculations are more involved than the
simple isotropic point-to-point interactions of PME, but are tractable and more
importantly \textit{local}, which contributes to scalability.  As we highlight
in the following sections, the short-range potential also has fast but algebraic
decay (up to $1/R^6$ in our examples), which is less attractive than the
exponential decay seen in PME, thus possibly leading to more local interactions.

Improvements to Ewald-type schemes range from coarsening strategies to reduce
the number of grid points by using a staggered mesh~\cite{Cerutti:2009} to
multilevel approaches~\cite{Cerutti:2010} that yield increased locality in the
FFT calculations while resulting in only a small increase in total work.
Moreover, other methods such as the Multilevel Summation
Method~(MSM)~\cite{Sagui:2001,Sandak:2001,Skeel:2002,Izaguirre:2005} take
different approach to operator splitting altogether.

In summary, the goal of this paper is to detail a method the incorporates
mesh-based screens and to investigate the accuracy of such an approach.  In
Section~\ref{sec:description}, we develop the mathematical construction of each
component.  In particular, we detail the screen functions that lead to the exact
sparse linear system for $\phsm$ and the local evaluation of $\phsr$.  In
Section~\ref{sec:cost}, we develop a performance model for the method and
discuss its implications in a parallel setting.  A numerical experiment is shown
in Section~\ref{sec:numerics} to confirm the accuracy of our method.  Additional
considerations and possible extensions are discussed in
Section~\ref{sec:discussion}.

\section{Description of method}\label{sec:description}

The Ewald decomposition is often viewed through the construction of a
screen potential to define the corresponding short-range and long-range
potentials.  The usual PME formulation is consistent with the original Ewald
decomposition in that it uses a Gaussian screen function
 \begin{equation}
  \rho^{\textnormal{G}}_{i}(\bx) = Q_i \rho_0 e^{-a^2 |\bx - \bx_i|^2}.
\end{equation}
This screen, as depicted in Figure~\ref{fig:splitting}, yields a short-range
potential so that $\phsr \propto \erfc(a|\bx-\bx_i|)/|\bx-\bx_i|$, which is straightforward to compute.
The resulting mesh potential satisfies the Poisson problem,
\begin{equation}
 - \nabla \cdot \nabla \phsm = \sum_i \rho_i(\bx),
 \label{eq:poisson}
\end{equation}
which is then optimally solved using FFTs on a mesh.

\begin{figure}[!ht]
 \begin{center}
 \includegraphics[trim=0.in 0.in 0.in 0.in, clip, width=0.8\textwidth]{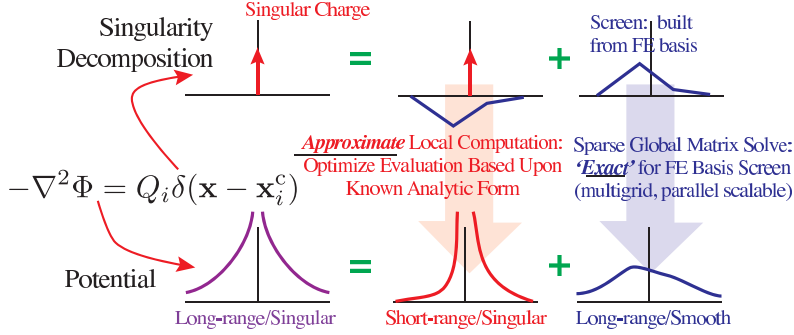}
  \caption{Introduction of a polynomial screen and resulting decomposition of the potential of the singular charge.}\label{fig:newsplitting}
\end{center}
\end{figure}

We instead propose screening functions $\rho_i(\bx)$ that are piecewise
polynomials of order $q$, as shown in Figure~\ref{fig:newsplitting}.
The corresponding potential is then solved exactly with (\ref{eq:poisson}) using
a finite element method with basis functions of order $p = q+2$.  That is, the
potential is represented exactly in the finite element space, making the
optimal resolution provided by an FFT-based solve unnecessary.

Next, we describe the details of the method, following the four
basic steps of \pppm~methods: assignment of charges to the mesh, solving for
the smooth potential on the mesh, transferring the potential back to the charge
locations, and calculating the point-to-point (short-range) interactions.   A
high-level synopsis of the algorithm is described in Algorithm~\ref{alg} to
illustrate the structural pieces of our approach.
\begin{algorithm}[!ht]
\SetKwComment{tcc}{\{}{\}}
\SetCommentSty{tiny}
\SetKwInput{Input}{Input} \SetKwInput{Output}{Return}
\SetKwIF{If}{ElseIf}{Else}{if}{}{else if}{else}{endif}
\SetKwFor{While}{while}{}{endw}
\SetKwFor{ForEach}{for each}{}{endfch}
\DontPrintSemicolon\SetVlineSkip{0.7in}
\Input{A mesh of elements $e_j$ and a group of point charges $Q_i$}
\Output{Potential at locations of charges}
\ForEach(\tcc*[f]{charge assignment, Section~\ref{sec:chargeassignment}}){$\mathrm{charge}\;Q_i$}{    place $Q_i$ in element\;
    solve for screen\tcc*[r]{(\ref{eq:momentmatrix}) or (\ref{eq:moment1d})}
    }
\ForEach{$\mathrm{element}\; e_j$}{       \If{$\mathrm{element}\; e_j \in \mathrm{surface}$}{     	   adjust boundary conditions as necessary\tcc*[r]{see (\ref{eq:bcadjust})}
       }
       apply charge assignment operator to form $\rho_\text{m}$\tcc*[r]{see (\ref{eq:defineassign})}
      }
   perform multigrid solve of $-\nabla^2 \Phi^\text{sm} = \rho_\text{m}$\tcc*[r]{see (\ref{eq:poisson})}
 \ForEach(\tcc*[f]{evaluations, Sections~\ref{sec:smootheval}~and~\ref{sec:shortrange}}){$ \mathrm{charge}\; Q_i$}{      $ \Phi_i \leftarrow \Phi_i^{\text{sr}} + \Phi_i^{\text{sm}}$\tcc*[r]{mesh-to-charge assignment}
     }
 \caption{Polynomial Screen Method for Calculating Potential}\label{alg}
\end{algorithm}

\begin{figure}[!ht]
  \begin{center}
    \includegraphics[trim=0.in 0.4in 0.in 0.in, clip, width=\textwidth]{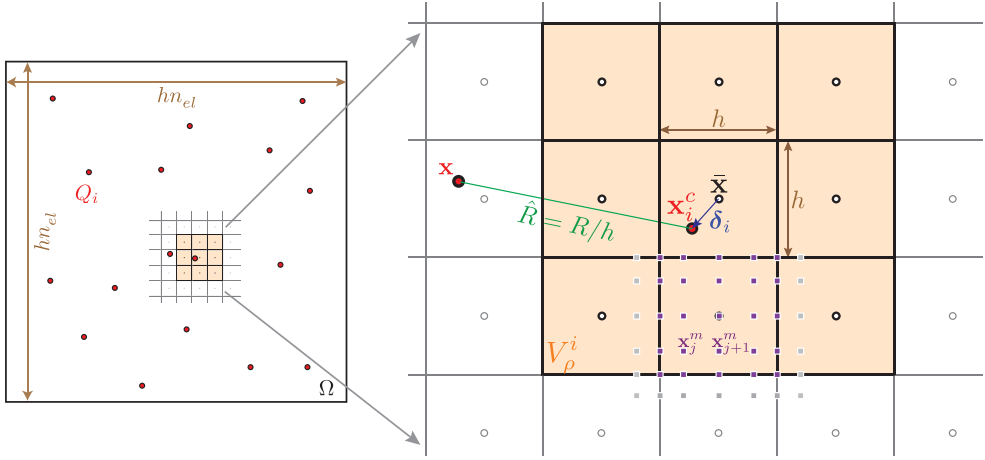}
    \caption{Schematic configuration showing $N$ charges of strength $Q_i$ at locations $\bx_i^\text{c}$ distributed in the cubic domain $\Omega$ of size $hn_\text{el} \times h n_\text{el} \times h n_\text{el}$, where $h$ is the size of the cubic finite elements.  The finite element centers are $\bar\bx$ and the collocation points are $\mathbf{x}_{j}^\text{m}$.}\label{fig:schematic}
  \end{center}
\end{figure}

We assume a collection of $N$ charges $\mQ \equiv {\{Q_i\}}_{i=1}^N$ located at
$\bxcharge_i$ in a cube $\Omega={[0,L]}^3$ (see Figure~\ref{fig:schematic}).  A
mesh with $\nel = \nelx\times\nely\times\nelz$ elements is
constructed to conform to the domain, and a uniform mesh is assumed in each
direction for simplicity of presentation~---~i.e., $n_\text{el} = \nelx = \nely = \nelz$.  Finally, the collocation points for $q$-order basis functions on
the mesh are denoted $\bxmesh_{j}$, with $j=1,\dots,M \equiv {(qn_\text{el} + 1)}^3$.

\subsection{Charge assignment}\label{sec:chargeassignment}

A central component of particle-mesh methods is the assignment of singular
charges to the mesh, yielding a mesh-based charge density function,
$\rhom(\bxmesh)$.  In particular, we seek an assignment function $\mW(\bx)$
that reflects our specially selected screen functions and provides a weighting
that distributes a charge $Q_i$ at $\bxcharge_i$ to each collocation point
$\bxmesh_j$ of the basis functions:
\begin{equation}
 \rhom(\bxmesh_j) =  Q_i\sum_{i = 1}^N \mW(\bxmesh_j;\bxcharge_i).
 \label{eq:assign}
\end{equation}
Existing methods use Lagrange  polynomials (PME~\cite{Darden:1993}) or
B-splines (\pppm~\cite{Hockney:1988} and SPME~\cite{Essmann:1995}) for this
weighting, the latter of which work particularly well with FFTs.  The charge
assignment function impacts both accuracy and efficiency of the method.  In our
approach we design an assignment operator based directly on polynomial
basis functions for compatibility with a finite-element-based Poisson solver.

\subsubsection{Defining the polynomial screens}

We define our screen density function for a single charge
$Q_i \in \mQ$ as
\begin{equation}
  \rho_i(\bx) = \sum_{j} c_j \scbasis_j(\bx),
\end{equation}
with linear superposition providing the extension to multiple charges.  Here
$\scbasis_j(\bx)$ are a collection of $q$-order Lagrange basis functions over an
index set determined as follows.  If charge $Q_i$ is located within element
$\tau_j$ of the mesh, we choose $V_{\rho}^i = \cup_{\tau \cap \tau_j \ne 0}
\tau$ to be the interpolation support of the charge assignment operator.  That
is, the support is the union of the element of the mesh that includes the charge
along with all neighboring elements, leading to a support of $27$ elements in
three dimensions.  Generalization to other choices for this support are briefly
discussed in Section~\ref{sec:discussion}.  To construct the polynomial screen,
we consider the degrees of freedom which are interior to or on the faces of the
element containing the charge.  For $q$-order interpolating polynomials, this
leads to $\dim(V_\rho^i) = {(q + 1)}^3$ degrees of freedom.  These degrees of
freedom are determined so that the charge-screen combination has a potential
that decays rapidly in space  by considering the multipole expansion of the
screen for a point well outside the screen, given by
\begin{equation}
\Phi^\text{sc}_i(\bx) =
  \frac{1}{\hat R}\int_{V_\rho^i} \rho_i(\boldsymbol{\xi})\;\mathrm{d}\boldsymbol{\xi}
- \frac{1}{\hat R^2}\int_{V_\rho^i}\rho_i(\boldsymbol{\xi})\;(\boldsymbol{\xi}\cdot\hat{\br})\; \mathrm{d}\boldsymbol{\xi}
+ \frac{1}{2\hat R^3}\int_{V_\rho^i} \rho_i(\boldsymbol{\xi})[3{(\boldsymbol{\xi}\cdot\hat{\br})}^2-|\boldsymbol{\xi}|^2] \;\mathrm{d}\boldsymbol{\xi}
+ \cdots
\label{eq:multipole}
\end{equation}
where $\hat{\bR} = (\bx - \bar{\bx}_\rho^i)/h$, with $\bx$ representing the
observation point, $\bar{\bx}_{\rho}^i$ is the center of the screen volume, and
$h$ is the mesh size.  The quantity $\hat{\br}$ is the unit direction
vector $\hat{\bR}/|\hat{\bR}|$.

For a charge $Q_i$ located at $\bx_i = (x_i,y_i,z_i)$ in element $\tau_j$ (see Figure~\ref{fig:schematic}),
we denote the offset $\bd_i = (\delta^x_i, \delta^y_i,
\delta^z_i) =  \bx_i- \bar{\bx}_j$ with respect to the center of the element $\bar{\bx}_j$, and
define the $(l,m,n)$-moment and \textit{centered} $(l,m,n)$-moment of the screen function as
\begin{align}
\rho_i^{(l,m,n)} &=  \int_{V_\rho^i} {(x-\delta^x_i)}^l
                                  {(y-\delta^y_i)}^m
                                  {(z-\delta^z_i)}^n\rho_i(\bx)\,\dx,\label{eq:moments}\\
\bar{\rho}_i^{(l,m,n)} &=  \int_{V_\rho^i} x^l y^m z^n \rho_i(\bx)\,\dx.\label{eq:cmoments}
\end{align}
where the origin is taken to be $\bar{\bx}_\rho^i$, the center of the screen volume $V_\rho^i$.
Dividing $\rho^{(0,0,0)}_i =
\int_{V^i_{\rho}} \rho_i(\bx)\,\dx $ by $\hat R = |\bx_i - \bx|/h$ gives the first term of the
screen's multipole expansion from (\ref{eq:multipole}). Thus, requiring $\rho_i^{(0,0,0)} = 1$ guarantees that the combined point-charge and screen have a potential that decays at least as fast as $1/R^2$ with distance from the point charge.  Likewise, zeroing higher moments of the screen enforces the cancellation of dipole and higher-order terms and further accelerates the long-range decay rate, thereby reducing the number of interactions that must be explicitly represented by point-to-point computations.  With the available degrees of freedom, a screen of order $q$ cancels all terms up
to $R^{-(q+1)}$, leaving $\phsr = 1/R -
\Phi^\text{sc} \sim R^{-(q+2)}$.  This is summarized in Table~\ref{tab:moments}, which shows the moments that result from performing the vector operations in the integrands of (\ref{eq:multipole}).
\begin{table}[ht!]
  \centering
\begin{tabular}{c  c  c }
\toprule
Power of $R$ & Single terms & Mixed terms \\
\midrule
$R^{-1}$ & $1$ & --- \\
$R^{-2}$ & $x, y, z$ & --- \\
$R^{-3}$ & $x^2,y^2,z^2$ & $xy, xz, yz$ \\
$R^{-4}$ & $x^3,y^3,z^3$ & $x^2y, xy^2, x^2z, xz^2, y^2z, yz^2, xyz$\\
$\vdots$ & $\vdots$ & $\vdots$ \\
$R^{-N}$ & $x^{N-1},y^{N-1},z^{N-1}$ &
\parbox[t][][t]{3.5cm}{\centering $x^l y^m z^n$,\;with $1\leq l,\,m,\,n \leq N-2$, and $l+m+n = N-1$}\\
\bottomrule
\end{tabular}
\caption{Polynomial terms in multipole expansion.}
\label{tab:moments}
\end{table}

\subsection*{Constructing the screen}

The goal is to perform multipole cancellations with screens
that are also compatible with the basis functions of our finite element
discretization.  With this description, each screen is composed of $\nsctotal =
\nsc^3 = {(q+1)}^3$ nodal screen basis functions, $\scbasis$:
\begin{equation}
   \rho_i(\bx) = \sum_{j = 0}^{\nsctotal-1} c_j\scbasis_j(\bx).
   \label{eq:sumscreen}
\end{equation}
Thus, revisiting (\ref{eq:assign}), the assignment operator $\mW$ is
\begin{equation}
 \mW(\bxmesh_j; \bxcharge_i)  =
  \begin{cases}
   \label{eq:defineassign}
       \sum_{k = 0}^{\nsctotal-1}\limits c_k(\bd_i)\scbasis_k(\bxmesh_j)
       &
       \qquad \bxmesh_j\in V_\rho^i \\
       \qquad 0 & \qquad \text{otherwise.}
  \end{cases}
\end{equation}
Restricting the first $\nsctotal$ moments leads to a $\nsctotal\times\nsctotal$ linear system
for the coefficients $\bc$ in (\ref{eq:sumscreen}):
\begin{equation}
  \renewcommand{\arraystretch}{1.5}
\begin{bmatrix}
\scbasis_0^{(0,0,0)} & \scbasis_1^{(0,0,0)} & \cdots & \scbasis_{\nsctotal-1}^{(0,0,0)} \\
\scbasis_0^{(1,0,0)} & \scbasis_1^{(1,0,0)}  & \cdots & \scbasis_{\nsctotal-1}^{(1,0,0)} \\
                                \vdots &  & \ddots & \vdots \\
\scbasis_0^{(q,q,q)} &  \scbasis_1^{(q,q,q)} &  \cdots & \scbasis_{\nsctotal-1}^{(q,q,q)} \\
\end{bmatrix}
\begin{bmatrix} c_0 \\
  c_1               \\
  \vdots            \\
  c_{\nsctotal-1}          \\
\end{bmatrix}
=
\begin{bmatrix}1 \\
                  0 \\
             \vdots \\
             0      \\
\end{bmatrix},
\label{eq:momentmatrix}
\end{equation}
where $\scbasis_j\,^{(l,m,n)}$ is the $(l,m,n)$-moment, as defined in
(\ref{eq:moments}) for $\rho$, of the $j$-th screen basis function
$\scbasis_j(\bx)$.

\begin{figure}[!ht]
  \psfrag{X}[][c]{$x$}
  \psfrag{Y}[c][1.0][1.0][-90]{$y$}
  \centering
  \begin{subfigure}[t]{0.45\textwidth}
    \includegraphics[trim=0.4in 0.4in 0.4in 0.4in, clip,width=0.8\textwidth]{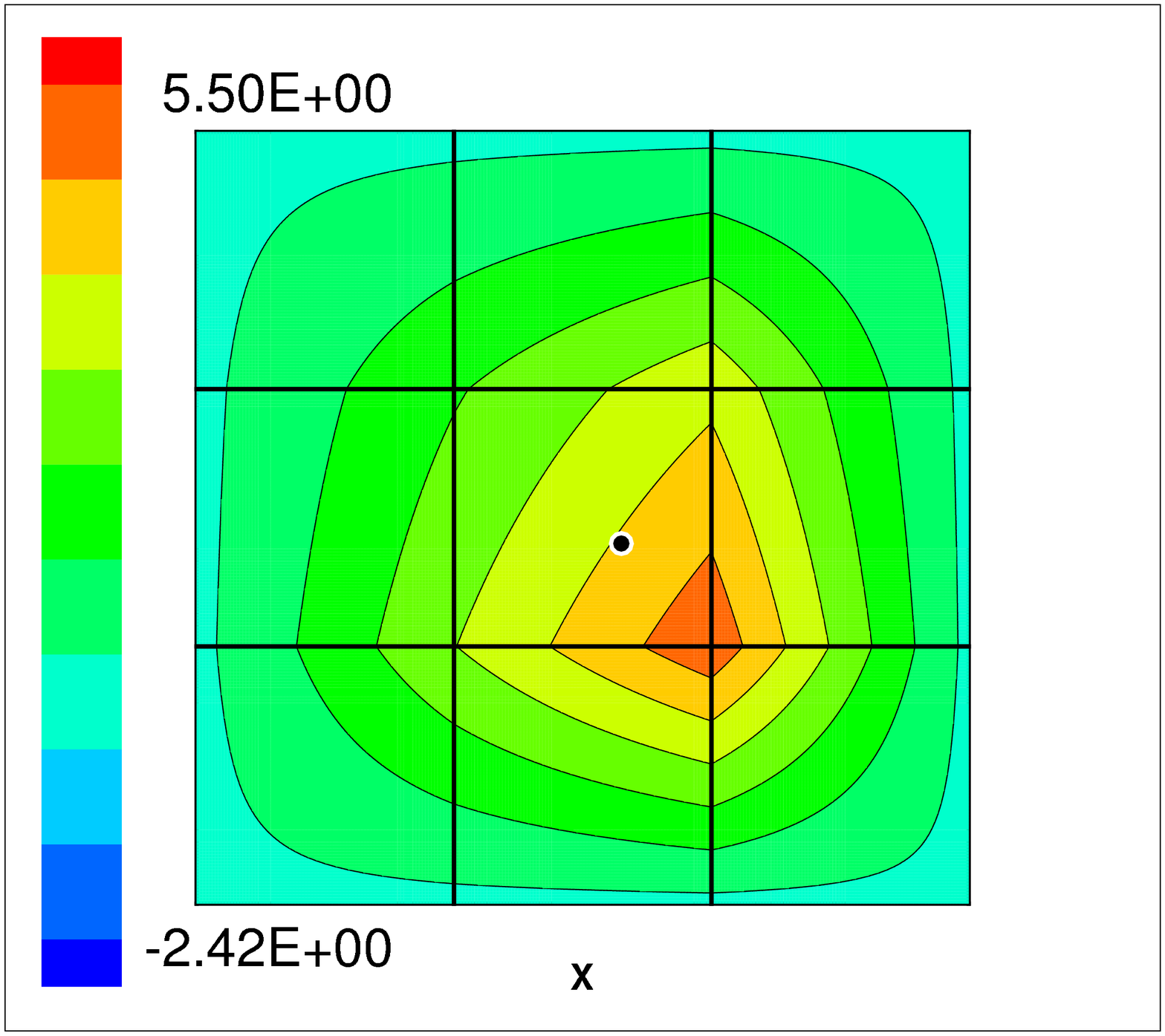}
    \caption{$q = 1$}
  \end{subfigure}
  \hfill
  \begin{subfigure}[t]{0.45\textwidth}
    \includegraphics[trim=0.4in 0.4in 0.4in 0.4in, clip,width=0.8\textwidth]{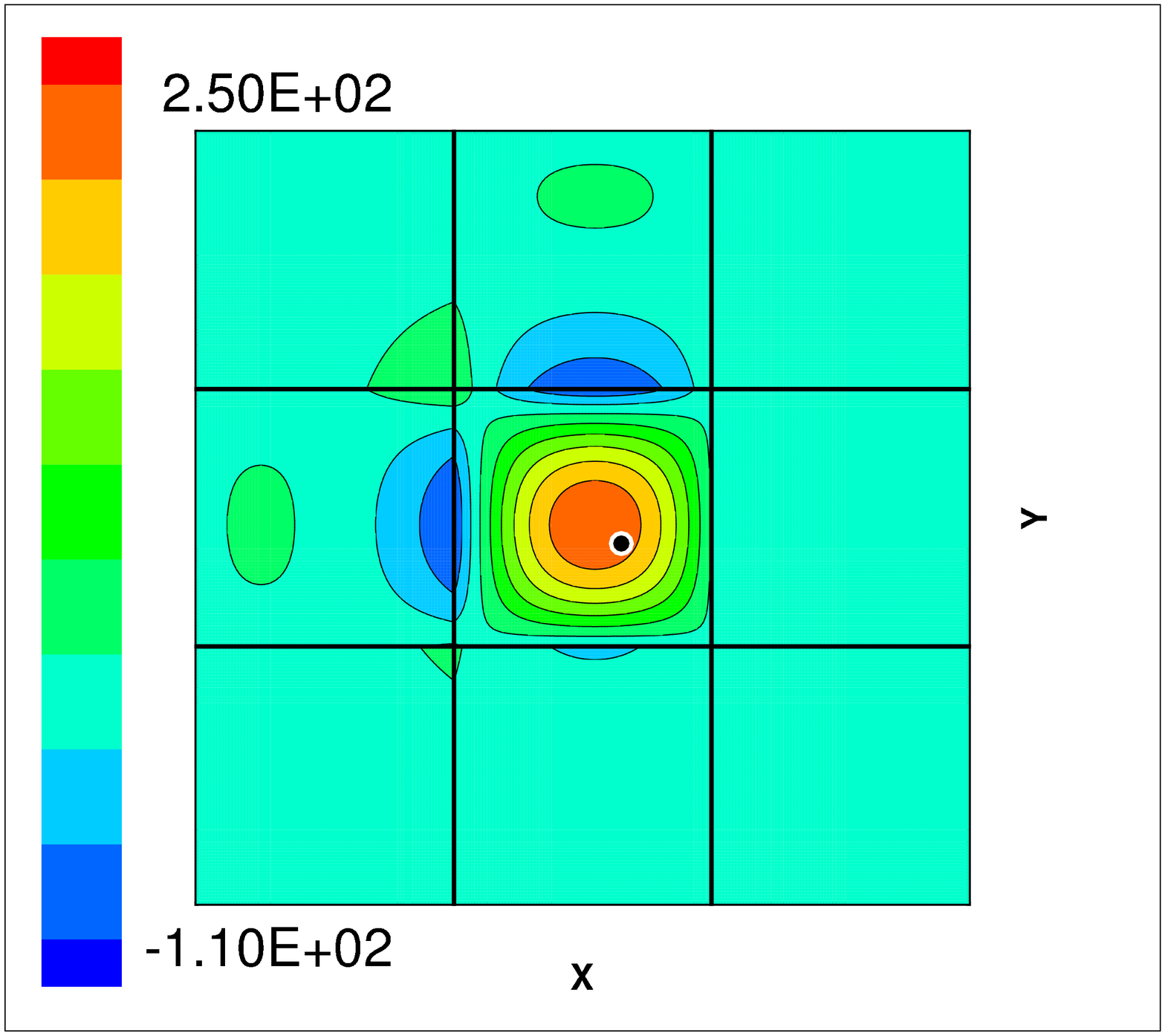}
    \caption{$q = 2$}
  \end{subfigure}
  \begin{subfigure}[t]{0.45\textwidth}
    \includegraphics[trim=0.4in 0.4in 0.4in 0.4in,  clip,width=0.8\textwidth]{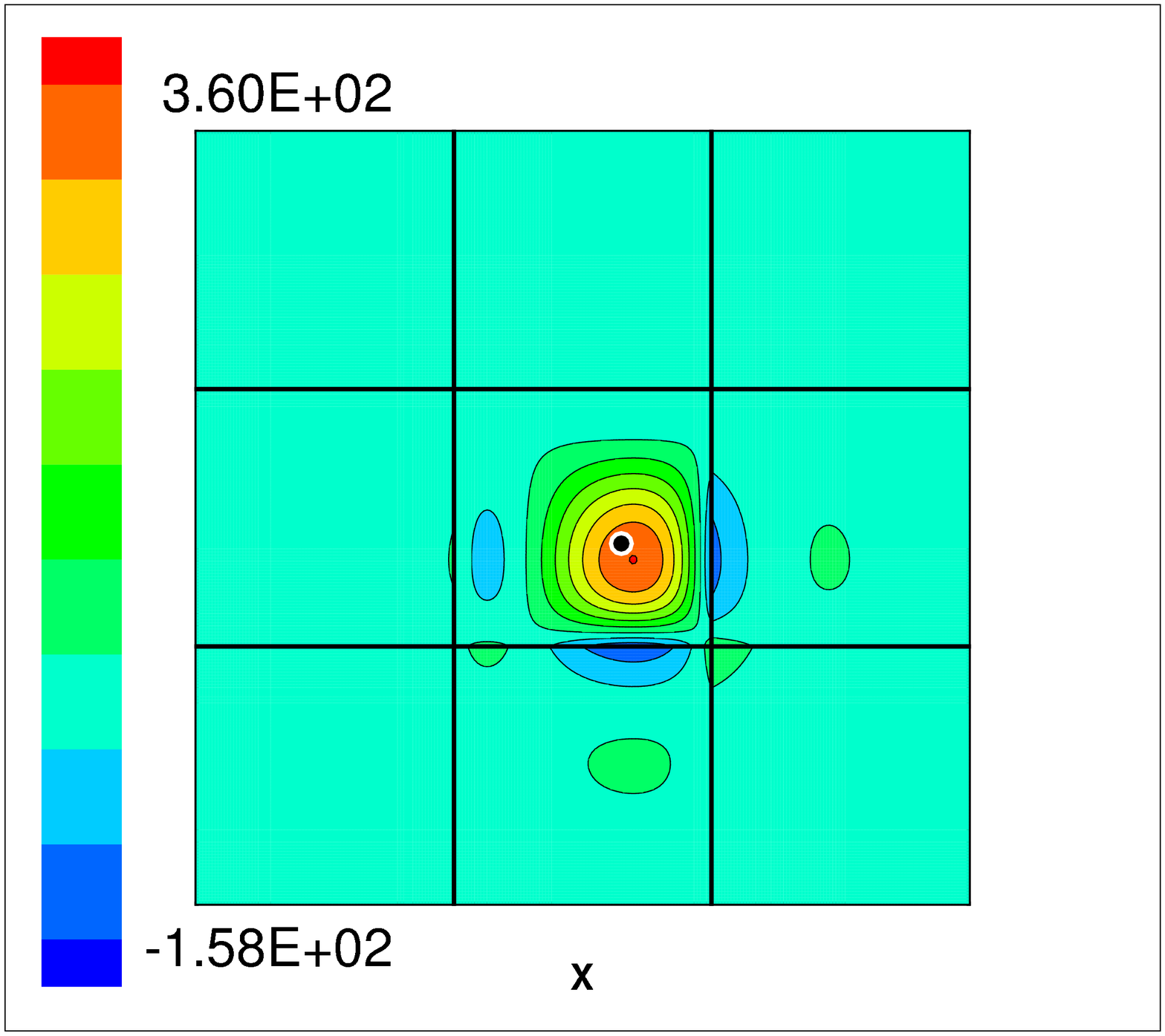}
    \caption{$q = 3$}
  \end{subfigure}
  \hfill
  \begin{subfigure}[t]{0.45\textwidth}
    \includegraphics[trim=0.4in 0.4in 0.4in 0.4in,  clip,width=0.8\textwidth]{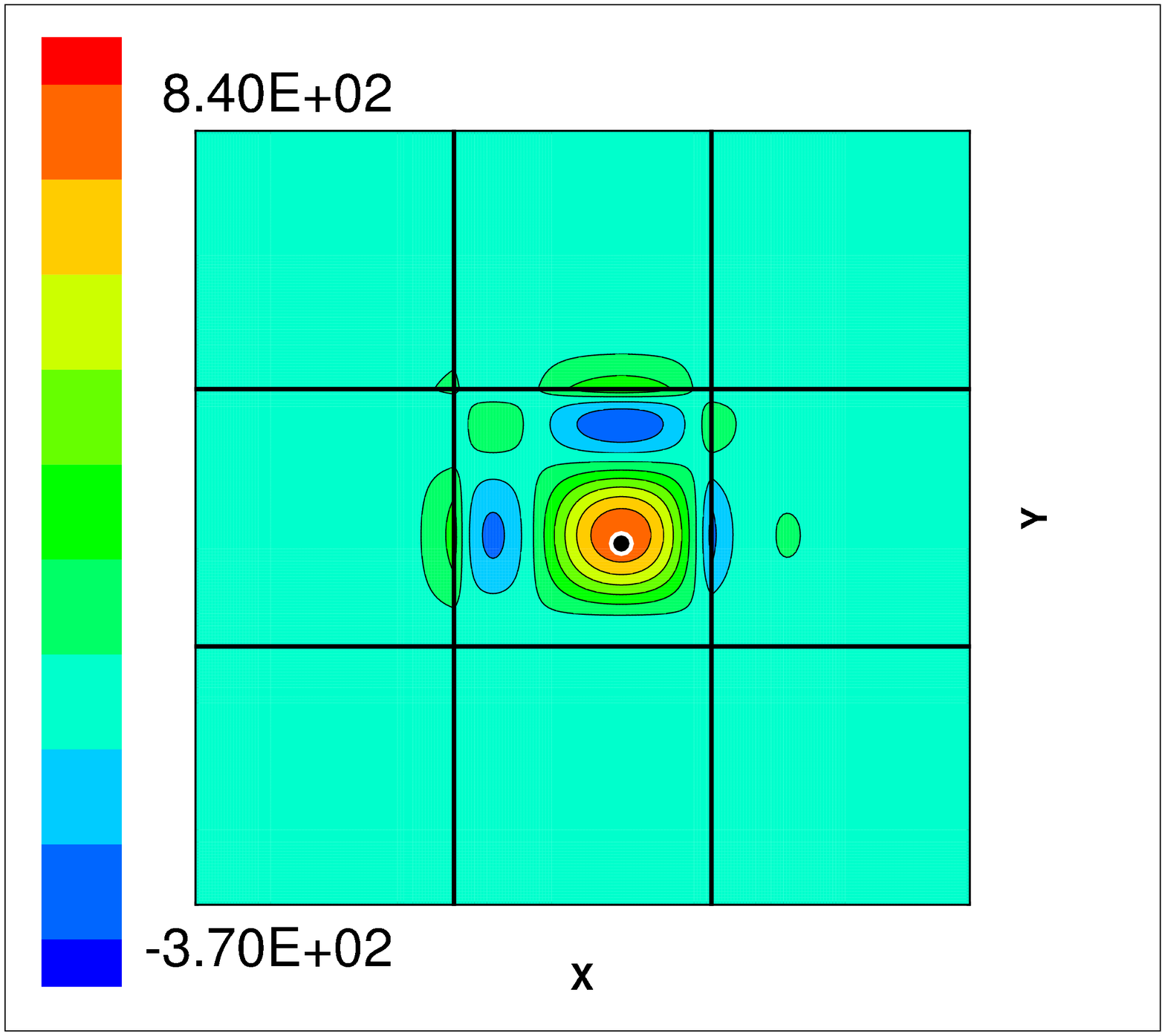}
    \caption{$q = 4$}
  \end{subfigure}
  \caption{Example linear through quartic screens for a single marked charge.}\label{fig:screens}
\end{figure}

Figure~\ref{fig:screens} shows cross-sections of example screens
constructed using $q = 1, \dots, 4$.  Each screen's peak is
attained near the marked charge, and the screens are constructed to decay to
zero at the edge of $V_\rho$.  The screens have support in the active screen region $V_\rho^i$ and for $q>1$, the screens are
in general non-monotone.

\begin{figure}[!ht]
\begin{center}
  \psfrag{R}[][c]{\raisebox{-0.2in}{$\hat R$}}
  \psfrag{L}[][c]{$\hat{R}$}
  \psfrag{S}[c][1.0][1.0][0]{$|\phsr|_\text{avg}$ or $|\phsr|_\text{max}$}
  \psfrag{P}[c][1.0][1.0][-90]{$\phsc$}
  \psfrag{a}[][l]{\hspace*{0.38in}\raisebox{-0.6in}{\begin{minipage}{0.15\textwidth}\scriptsize$q = 1$\\$q = 2$\\$q = 3$\\$q = 4$\end{minipage}}}
  \psfrag{b}[][c]{}
  \psfrag{c}[][c]{}
  \psfrag{d}[][c]{}
   \psfrag{z}[][c]{$\quad1/R$}
   \begin{subfigure}[t]{0.45\textwidth}
         \includegraphics[trim=0.2in 0.2in 0.2in 0.2in, clip, width=\textwidth]{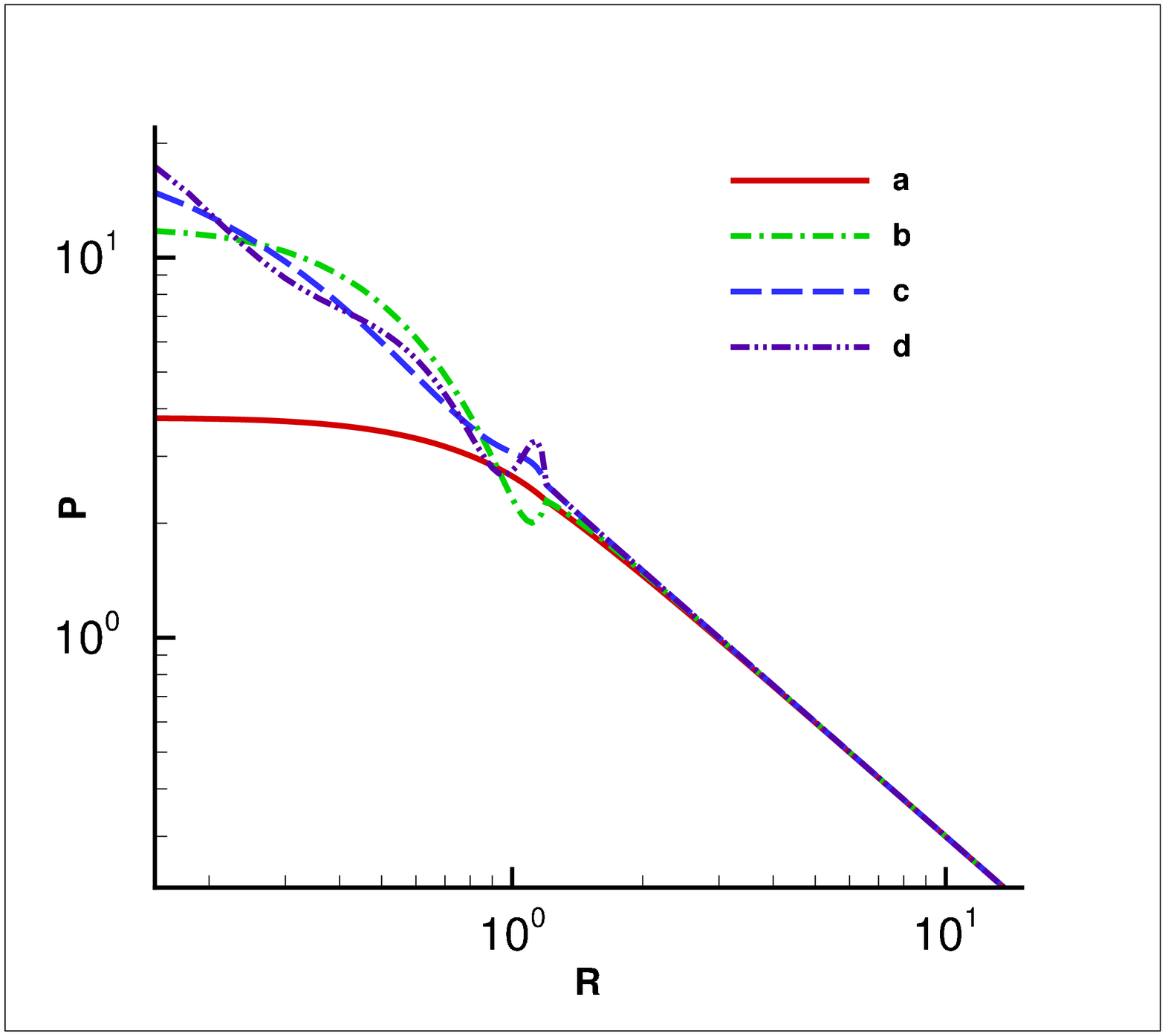}
         \caption{}
     \end{subfigure}  \hfill
    \begin{subfigure}[t]{0.45\textwidth}
       \includegraphics[trim=0.2in 0.2in 0.2in 0.2in, clip, width=\textwidth]{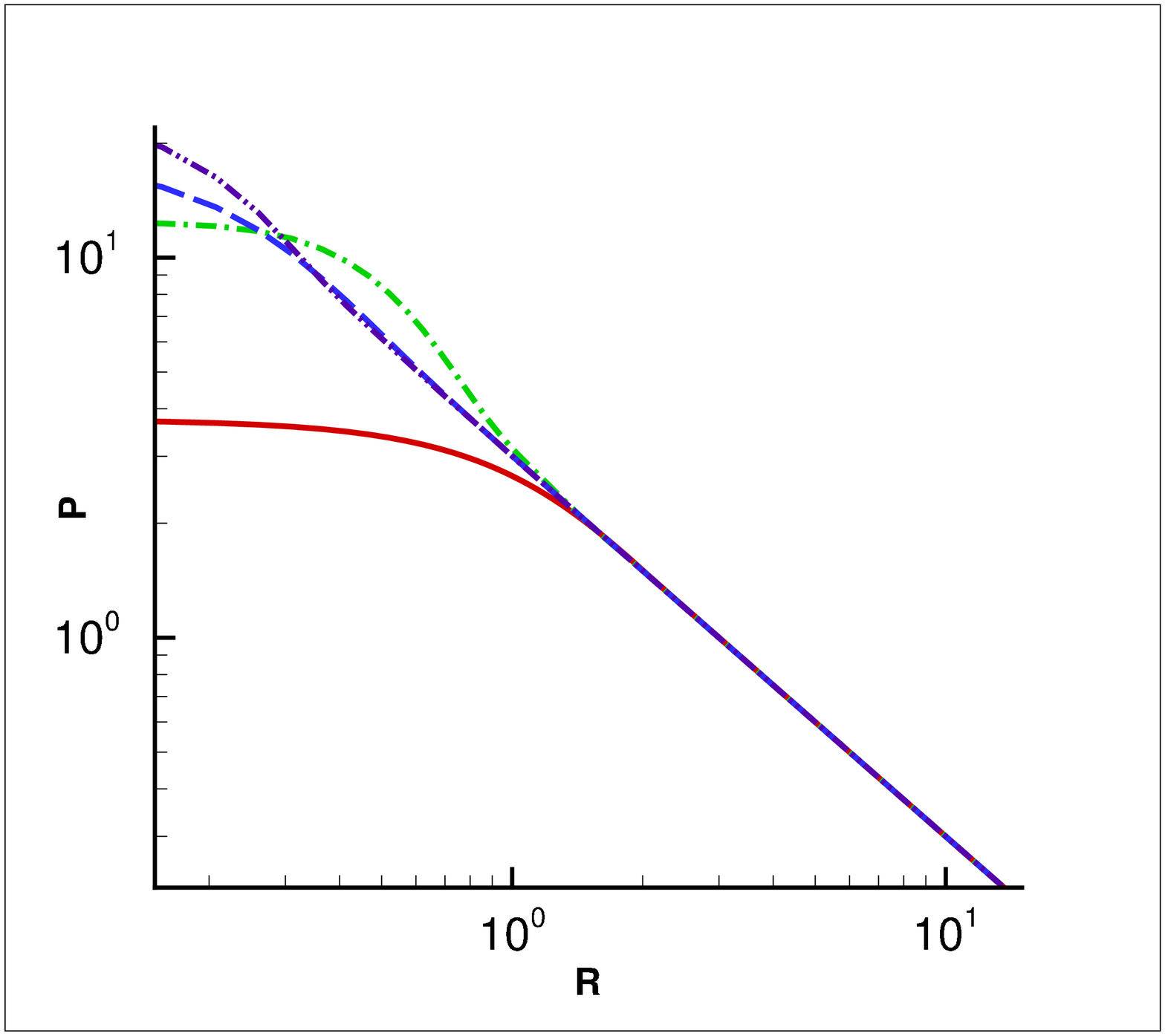}
       \caption{}
    \end{subfigure}
    \caption{Screen potential in two directions: (a) $\bx-\bx_i \propto (1, 0, 0)$ and (b) $\bx-\bx_i \propto (1, 1, 1)$.}\label{fig:screenpo}
\end{center}
\end{figure}

\begin{figure}[!ht]
 \psfrag{R}[][c]{$R$}
  \psfrag{L}[][c]{$\hat{R}$}
  \psfrag{S}[c][1.0][1.0][0]{\begin{minipage}{0.3\textwidth}\centering$|\phsr|_\text{avg}$ or $|\phsr|_\text{max}$\\ $\quad$ \end{minipage}}
  \psfrag{P}[c][1.0][1.0][-90]{$\phsc$}
    \psfrag{a}[][r]{\hspace*{1.1in}\raisebox{-0.25in}{\scriptsize$q = 1,\; \phsr \sim \hat{R}^{-3}$}}
  \psfrag{b}[][]{\hspace*{1.1in}\raisebox{-0.2in}{\scriptsize$q = 2,\; \phsr \sim \hat{R}^{-4}$}}
  \psfrag{c}[][]{\hspace*{1.1in}\raisebox{-0.25in}{\scriptsize$q = 3,\; \phsr \sim \hat{R}^{-5}$}}
  \psfrag{d}[][]{\hspace*{1.1in}\raisebox{-0.35in}{\scriptsize$q = 4,\; \phsr \sim \hat{R}^{-6}$}}

 \begin{center}
    \includegraphics[trim=0.2in 0.2in 0.2in 0.2in, clip, width=0.6\textwidth]{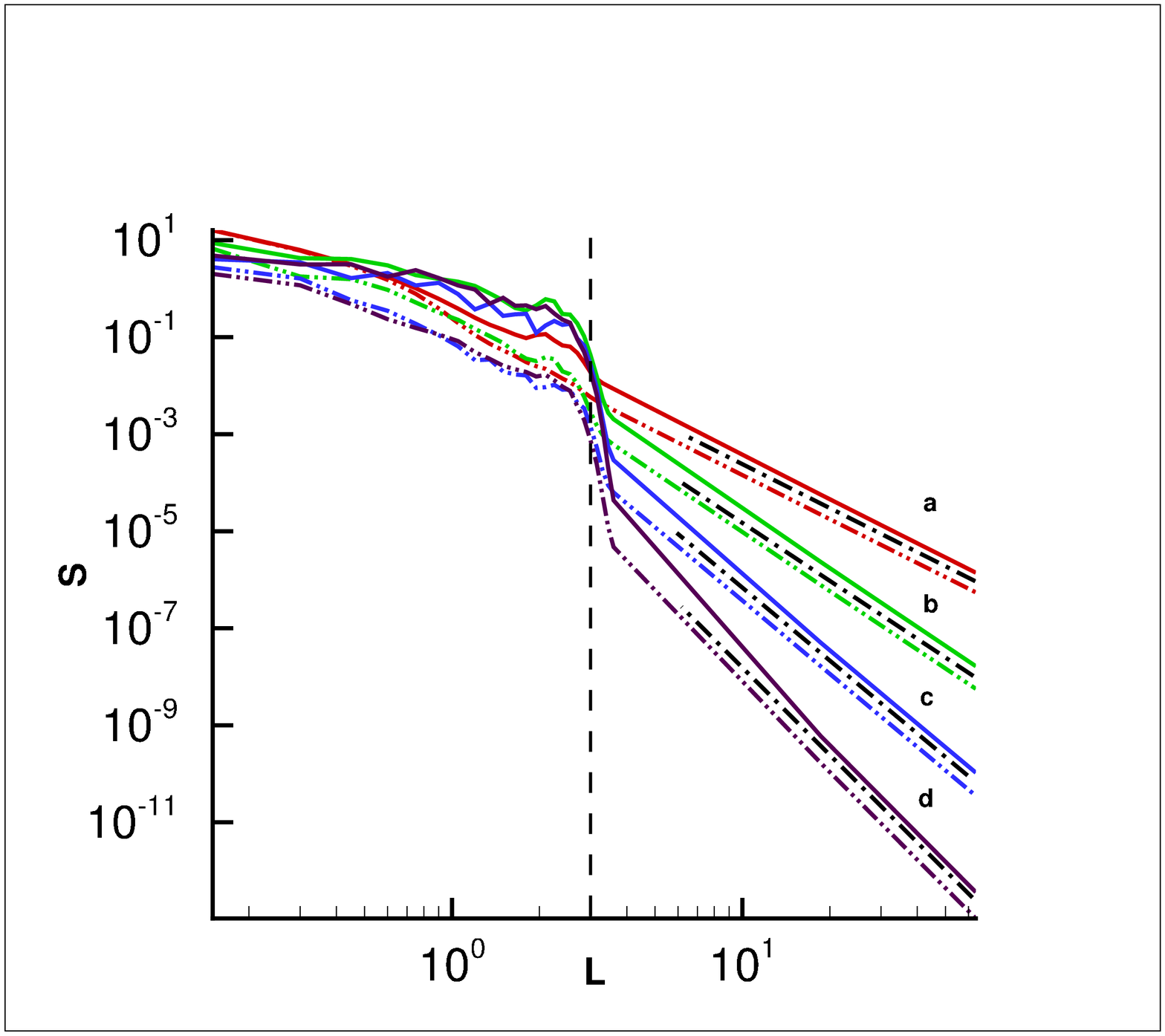}
   \caption{Maximum (solid) and average (dash-dot-dot) short-range potentials for screens with $q = 1$ to $4$. The straight dash-dot lines denote the expected slopes. }\label{fig:decayavg}
\end{center}
\end{figure}

We confirm in Figure~\ref{fig:screenpo} that screens constructed in this
fashion  yield potentials with the expected behavior: in all cases, the
far-field behavior approximates $1/R$.   In Figure~\ref{fig:decayavg}, we
estimate mean and peak errors incurred for point-to-point interaction
truncation at a distance $\hat R_c = R_c/h$.  To do this, the short-range
potential is constructed for $N=84$ charge locations and sampled in 42
directions; the behavior of $|\phsr|$ is shown in the figure as the weighted
average $|\phsr|_\text{avg}$ and the maximum from all samples
$|\phsr|_\text{max}$.  In addition, the radius from a charge location is
also normalized as $\hat{R} = R/h$, and is denoted $\hat{R}$ (see
Figure~\ref{fig:schematic}).  We see that twice the screen-size scale $\hat{R}
\approx 3$ corresponds to the expected start of the asymptotic decay behavior.
This is the distance at which a multipole expansion is generally considered
``well-separated'' and expected to show convergence with $R$.  Both the mean
and peak errors show the expected behavior for increasing $q$ beyond this
distance.

\subsection*{A fast algorithm for screen construction}

Solving (\ref{eq:momentmatrix}) directly requires $\mO(\nsctotal^3)$ operations for
each screen, which is feasible, but is not necessary in general. In the following, we
design a fast algorithm for computing screens, which follows from a
generalization of the Parallel Axis Theorem applied to moments used in the
system.  First, we note that the moments are additive.  For example, the first moment in
variable $x$ of basis function $\scbasis_j$ satisfies
\begin{equation}
\begin{split}
  \scbasis_j^{(1,0,0)}
  &= \int_{V_\rho} x\scbasis(\bx)\;\dx - \delta^x_i \int_{V_\rho}\scbasis(\bx)\;\dx\\
  &= \bscbasis_{j}^{(1,0,0)} - \delta^x_i\scbasis_j^{(0,0,0)},
\end{split}
\end{equation}
where $\bscbasis_{j}\,^{(1,0,0)}$ is the \textit{centered} moment of the
$j$-th basis function as in (\ref{eq:cmoments}).
Therefore, second row of (\ref{eq:momentmatrix}) is equivalent to
\begin{equation}
\sum_{j = 0}^{\nsctotal-1}c_j\bscbasis_j^{(1,0,0)}
= \delta^x_i\sum_{j = 0}^{\nsctotal-1} c_j\scbasis_j^{(0,0,0)}
= \delta^x_i.
\end{equation}
Similarly, the second moment of $\scbasis_j$ satisfies
\begin{equation}
\begin{split}
  \scbasis_j^{(2,0,0)}
  &= \int_{V_\rho} x^2\scbasis(\bx)\;\dx - 2\delta^x_i \int_{V_\rho}x\scbasis(\bx)\;\dx + {(\delta^x_i)}^2\int_{V_\rho} \scbasis(\bx)\; \dx \\
  &= \bscbasis_{j}^{(2,0,0)} - 2\delta^x_i \bscbasis_{j}^{(1,0,0)} + {(\delta^x_i)}^2\scbasis_j^{(0,0,0)},
  \end{split}
\end{equation}
so the third row of (\ref{eq:momentmatrix}) becomes
\begin{equation}
\begin{split}
\sum_{j = 0}^{\nsctotal-1} c_j \bscbasis_{j}^{(2,0,0)}
& = 2\delta^x_i\sum_{j=0}^{\nsctotal-1} c_j\bscbasis_{j}^{(1,0,0)} - {(\delta^x_i)}^2\sum_{j = 0}^{\nsctotal-1} c_j \scbasis_j^{(0,0,0)}\\
& = 2 {(\delta^x_i)}^2 - {(\delta^x_i)}^2= {(\delta^x_i)}^2.
\end{split}
\end{equation}
Continuing this procedure for other moments in (\ref{eq:momentmatrix}) yields:
\begin{equation}
  \renewcommand{\arraystretch}{1.8}
\begin{bmatrix}
    \scbasis_{0}^{(0,0,0)} & \scbasis_1^{(0,0,0)} & \hdots  & $\quad$ & \scbasis_{\nsctotal-1}^{(0,0,0)} \\
    \bscbasis_{0}^{(1,0,0)} & \bscbasis_{1}^{(1,0,0)}&           & $\quad$   &                    \\
  \vdots                                &                                        &              & $\quad$   &      \vdots           \\
  \bscbasis_{0}^{(q,0,0)} &                                        &   \ddots  & $\quad$  & \bscbasis_{{\nsctotal-1}}^{(q,0,0)}\\
  \bscbasis_{0}^{(0,1,0)} &                                        &                & $\quad$   &                    \\
  \vdots                               &                                        &                 & $\quad$  &     \vdots        \\
  \bscbasis_{0}^{(q,q,q)} &                                         &  \hdots        & $\quad$  & \bscbasis_{{\nsctotal-1}}^{(q,q,q)} \\
\end{bmatrix}
\begin{bmatrix}
   c_0 \\
   c_1 \\
        \\
   \vdots    \\
   \\
        \\
   c_{\nsctotal-1} \\
\end{bmatrix}
=
\begin{bmatrix}
   1 \\
   \delta^x_i \\
    \vdots \\
   {(\delta^x_i)}^q \\
    \delta^y_i   \\
      \vdots  \\
    {(\delta^x_i)}^q{(\delta^y_i)}^q{(\delta^z_i)}^q   \\
\end{bmatrix},
\label{eq:momentmatrixcentered}
\end{equation}
which we write compactly as $\mC\bc = \bbf$.
An advantage of this form is that for a uniform mesh, the matrix $\mC$ is the same for each
screen, since the moments reference the center of the element.  As a result, the
matrix is pre-factorized leading to a complexity of only $\mO(\nsctotal^2)$ to solve for each
screen.

For small $q$ this yields a small operation count, yet the computation is further reduced
if $\psi(\bx)$ is separable, as is the case for the regular cubic mesh shown in Figure~\ref{fig:schematic}.  In this case,
\begin{equation}
 \scbasis_\kappa(\bx) =  \omega_{i}(x)\omega_{j}(y)\omega_{k}(z),
\label{eq:oned}
\end{equation}
where $\omega_{i}$ are the one-dimensional nodal basis functions for a mesh size $h$ and $\kappa = i + (q+1)j + {(q+1)}^2k$ with $i,j,k \in [0,q]$.
With $\psi$ separable, the moment integrals are also separable:
\begin{align}
\bar{\scbasis}_\kappa^{(l,m,n)} =&  \int_{V_\rho^i} x^l y^m z^n \scbasis_\kappa(\bx)\,\dx \nonumber \\
=&\left(\int_{-3h/2}^{3h/2} x^l \omega_{i}(x) \, \mathrm{d}x \right)\left(\int_{-3h/2}^{3h/2} y^m \omega_{j}(y)  \, \mathrm{d}y\right)\left( \int_{-3h/2}^{3h/2}z^n \omega_{k}(z)\,\mathrm{d}z \right).
\label{eq:cmomentsbasis}
\end{align}
Following the notation of (\ref{eq:cmoments}), we define
\begin{equation}
\bar{\omega}_{i}^{(l)} = \int_{-3h/2}^{3h/2} x^l \omega_{i}(x) \, \mathrm{d}x
\end{equation}
and likewise for the one-dimensional $y$ and $z$ centered moments.
We take $c_\kappa = w^x_i w^y_j w^z_k$ and recognize that the right-hand side of (\ref{eq:momentmatrixcentered}) is also separable as  $f_\mu = {(\delta^x)}^l {(\delta^y)}^m {(\delta^z)}^n$ with $\mu =  l + (q+1)m + {(q+1)}^2n$.  This yields three equivalent ${(q+1)}^2$ systems of the form
\begin{equation}
\sum_{i = 0}^q \bar{\omega}_{i}^{(l)}(x)w_i^x = {(\delta^x)}^l, \quad l = 0,\hdots,q
\label{eq:moment1d}
\end{equation}
 which are solved independently.  The inverse of the matrix is computed once and applied for all right hand sides, resulting in only $\mO(\nsc^2)$ operations per screen.  This method also extends to regular rectangular meshes, where $h$ is not necessarily equal in each direction.

\subsection{Solution of the mesh potential} \label{sec:meshsolve}

Given our  construction of the screen $\rho$ using the finite element basis functions via (\ref{eq:defineassign}), the solution of $\Phi^\text{sm}$ is straightforward.  We simply use a finite element
solver with basis functions $\mbasis$ of
order $p = q + 2$, which results in a symmetric, positive definite sparse linear
system (under realistic assumptions regarding the boundary conditions) that does not introduce any numerical approximations.

\begin{figure}[!ht]
  \centering
  \psfrag{Y}[c][1.0][1.0][-0]{$\|r\|_2$}
  \psfrag{x}[c][1.0][1.0][0]{\raisebox{-0.in}{CG iterations}}
  \psfrag{a}[][l]{\raisebox{-0.42in}{\hspace*{1.in}\begin{minipage}{0.4\textwidth}
        conjugate gradient (CG) \\
        multigrid preconditioned CG \end{minipage}}}
  \psfrag{b}[][]{increasing $p$}
  \includegraphics[trim=0.2in 2.4in 0.2in 0.2in,clip,width=0.85\textwidth]{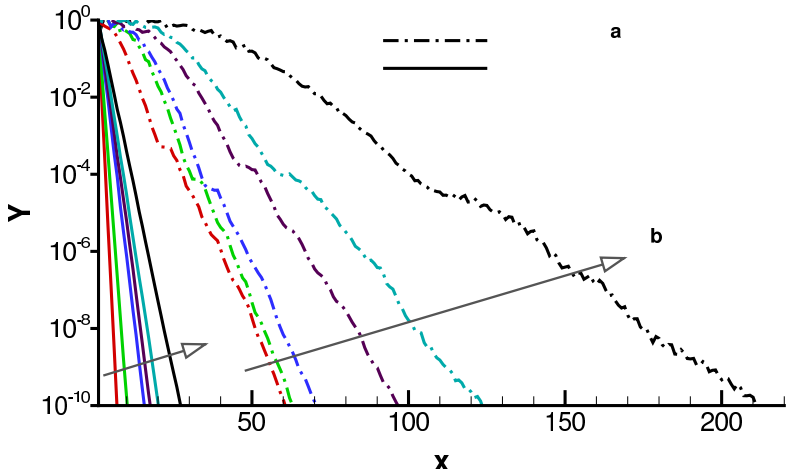}
  \caption{Convergence history for the Poisson problem using basis elements of order $p=1,\dots,6$ for both the conjugate gradient method (dot-dashed) and multigrid preconditioning (solid).}\label{fig:mghighorder}
\end{figure}

Multigrid preconditioners are effective for this problem, even for high-order
bases, and allow the sparse matrix problem to be solved to any level of
accuracy.  As an example, consider the case of a high-order finite element
discretization of the Poisson problem with Dirichlet boundary conditions.
Figure~\ref{fig:mghighorder} shows the convergence history of a multigrid
preconditioned conjugate gradient method for basis functions of order $p = 1$
through $6$.  An algebraic multigrid preconditioner, based on smoothed
aggregation using a more general strength measure~\cite{2011_OlScTu_energymin}
and optimal interpolation operator~\cite{2010_OlScTu_evosoc}, is used.  We
observe only a weak dependence on $p$.  Moreover, more advanced multigrid
techniques have shown still better scalings for both Poisson and other elliptic
problems such as for Stokes flow~\cite{2005_HeMaMcOl_hoamg,
2007_Olson_hoamgtri}.  Importantly for our principal objective, multigrid
preconditioners are well-known to exhibit high parallel
efficiency~\cite{2010_GhGr_feasibility, 2012_BaAl_etal_Scalability}.  In the
following tests, we use AMG through the BoomerAMG package~\cite{boomeramg}.

\subsection{Evaluation of the smooth potential at charge locations} \label{sec:smootheval}

The next step is to evaluate $\phsm$
at the charge locations.  For standard \pppm~and PME implementations, this involves
interpolation with the Lagrangian or B-spline basis functions from the charge
assignment.  In contrast, our method requires no interpolation, though interpolation can be used to speed up calculations, if desired.
Since the smooth potential exists in each element as a linear
combination of coefficients~---~i.e., the values of $\phsm$ at $\bxmesh$ for all
$\nmesh$ points in an element~---~and the basis functions $\mbasis$,
the smooth potential is expressed \textit{exactly} at any point as
\begin{equation}
 \phsm(\bx) = \sum_{j = 1}^{\nmesh}  \phsm(\bx^{\text{m}}_j) \mbasis_j(\bx-\bar{\bx}),
  \label{eq:evalsmooth}
\end{equation}
where $\nmesh$ is the number of collocation points in an element. Direct evaluation at the charge locations $\bx = \bx_i^\text{c}$ is straightforward.

\subsection{Short-range potential}   \label{sec:shortrange}

Our formulation for the exact mesh solution yields a more complex short-range
interaction than PME.\@  In addition to $R$, the short-range interaction now
also depends on the position of the charge relative to the underlying mesh.
Consequently, additional effort is required to evaluate the short-range
interaction.  However, the calculation is \textit{local}, so it does not
inhibit parallel efficiency.

The short-range potential at point $\bx$ due to a charge $Q_i$
located at $\bx_i$ is
\begin{equation}
   \phsr_i(\bx) = \frac{Q_i}{|\bx-\bx_i|} - \phsc_i(\bx),
\end{equation}
where
\begin{equation}
    \phsc_i(\bx) = \int_{V_\rho^i} \frac{\rho_i(\bxi)}{| \bx - \bxi | } \mathrm{d}\bxi.
    \label{eq:screenintegral}
\end{equation}
Though feasible, performing accurate quadrature for each screen individually
is computationally expensive.  We therefore shift a significant portion of this
computational effort to a pre-processing step, for which there are multiple options.

One approach is to consider a look-up table of
pre-computed values for the screen potential evaluated at $\bx_j$ due
to a charge at $\bx^\text{c}_i$.  These values are represented in a six-dimensional look-up table as $\phsc(\bx_j-\bxcharge_i; \bd_i)$, since
they are a function of the difference between the evaluation point and the
charge location, and also the offset of the charge within its element (which
determines the screen).

With some additional computation, but still without resorting to direct evaluation of (\ref{eq:screenintegral}), it is possible to remove the charge offset interpolation to reduce errors.   We accomplish this by recognizing the screen's formulation as a linear combination of basis functions,
\begin{align}
\phsc_i(\bx) &= Q_i\sum_{j= 0}^{\nsctotal-1}c_j(\bd_i)
                \int_{V_\rho^i} \frac{\scbasis_j(\bxi)}{| \bx-\bxi |} \mathrm{d}\bxi \\
             &=  Q_i\sum_{j= 0}^{\nsctotal-1} c_j(\bd_i) \Phi_j^{\text{b,sc}}(\bx-\bx^\text{c}_i).
  \label{eq:sc}
\end{align}
This approach yields $\nsctotal$ look-up tables for
basis-function potential values $\Phi^{\text{b,sc}}(\bx_j-\bx^\text{c}_i)$.  However, for $q \ge 2$ the
polynomial nature of the screen leads to non-monotonic decay for some directions within the region
where the screen is active, as shown in
Figure~\ref{fig:screenpo}.  Consequently, a direct implementation
of a look-up table for such functions requires sufficient resolution, which is harder to achieve for larger $q$.
For good performance, knowledge of the underlying structure of the screen potentials
should be used to inform both
the storage locations for the look-up table values and the interpolation method.

\subsection{A note about the self term}
If the point $\bx$ is the location of a charge, we do not wish to include the potential due to this charge in our calculation.  However, we do still need to subtract the screen potential from the charge's own screen, which is sometimes called the ``self'' term.   We can allow for this by amending our short-range potential expression to include both cases:
\begin{equation}
  \phsr_i(\bx) =
\begin{cases}
 \frac{Q_i}{|\bx-\bx_i|} - \phsc_i(\bx) & \qquad |\bx-\bx_i| > 0 \\
                             - \phsc_i(\bx) & \qquad \text{otherwise.}
  \end{cases}
  \end{equation}

\subsection{Alternate boundary conditions}
\label{sec:bc}

The formulation above is presented under the assumption of periodic boundary
conditions, which is the simplest case and important for a range of
applications.  It is straightforward to generalize boundary conditions via the
mesh potential $\phsm$.  This is accomplished by adjusting for short-range
effects present at the boundary and then proceeding in the usual manner for a
finite element problem with the given type of boundary conditions.  For
example, for a Dirichlet boundary condition of $\Phi = g$ on $\partial V$, the
condition for our mesh problem becomes
\begin{equation}
  \phsm\rvert_{\partial V} = g - \phsr \rvert_{\partial V},
  \label{eq:bcadjust}
\end{equation}
which leads to
\begin{equation}
  \Phi \rvert_{\partial V} = \phsm\rvert_{\partial V} + \phsr  \rvert_{\partial V} = g.
\end{equation}
A similar approach is used in~\cite{Hernandez-Ortiz:2007}.
 This also extends to the case of
Neumann or mixed-type boundary conditions, with the usual constraint to address the non-uniqueness of the fully Neumann problem.  Free-space conditions impose the
usual challenges but are no more difficult for the proposed scheme than for any mesh-based
Poisson solver.

\section{Performance model} \label{sec:cost}

\begin{figure}[h]
 \psfrag{T}[l][1.0][1.0][-90]{$\!\!\!\!\!\!\!\!\!\!\!\!\!\!\!\!t$ (sec)}
  \psfrag{N}[c][1.0][1.0][0]{\raisebox{-0.2in}{$N$}}
  \psfrag{p}[c][1.0][1.0][0]{\raisebox{-0.1in}{$q$}}
   \psfrag{M}[c][1.0][1.0][0]{\raisebox{-0.2in}{$M$}}
   \psfrag{a}[][l]{\hspace*{0.58in}\raisebox{-0.64in}{\begin{minipage}{0.25\textwidth}\scriptsize$q = 4$\vspace*{0.02in}\\$q = 3$\vspace*{0.02in}\\$q = 2$\vspace*{0.02in}\\$q = 1$\end{minipage}}}

   \begin{subfigure}[t]{0.48\textwidth}
    \includegraphics[trim=0.5in 0.2in 0.5in 0.2in, clip, width=\textwidth]{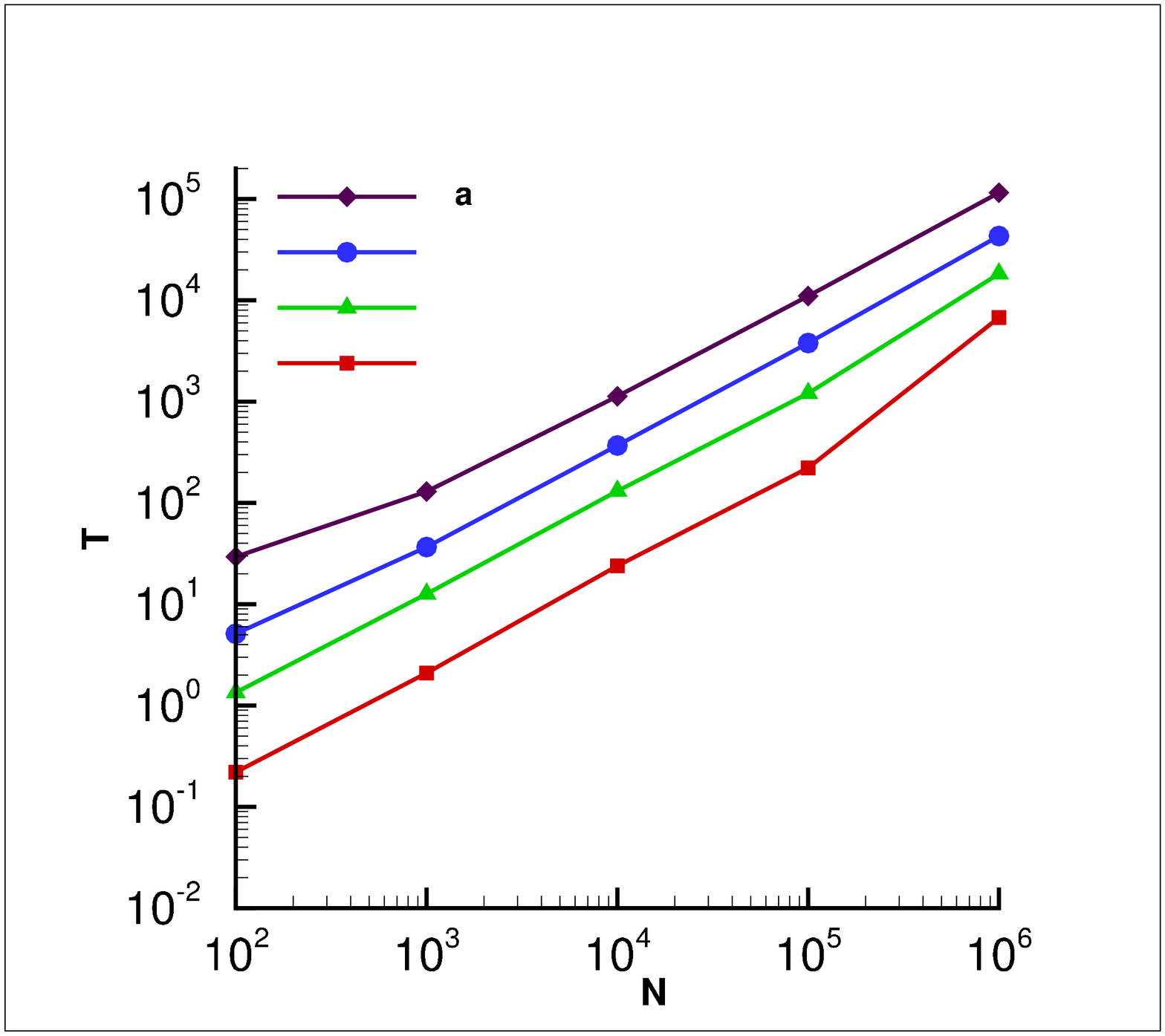}

  \caption{Total CPU time vs.\ $N$ for the main $N$-dependent components of the algorithm, including screen creation, short-range calculation, basis function evaluation, and combination of short-range and mesh potentials at charge locations.}\label{fig:timingstotalN}
  \end{subfigure}\hfill
\begin{subfigure}[t]{0.48\textwidth}
    \includegraphics[trim=0.2in 0.2in 0.6in 0.7in, clip, width=\textwidth]{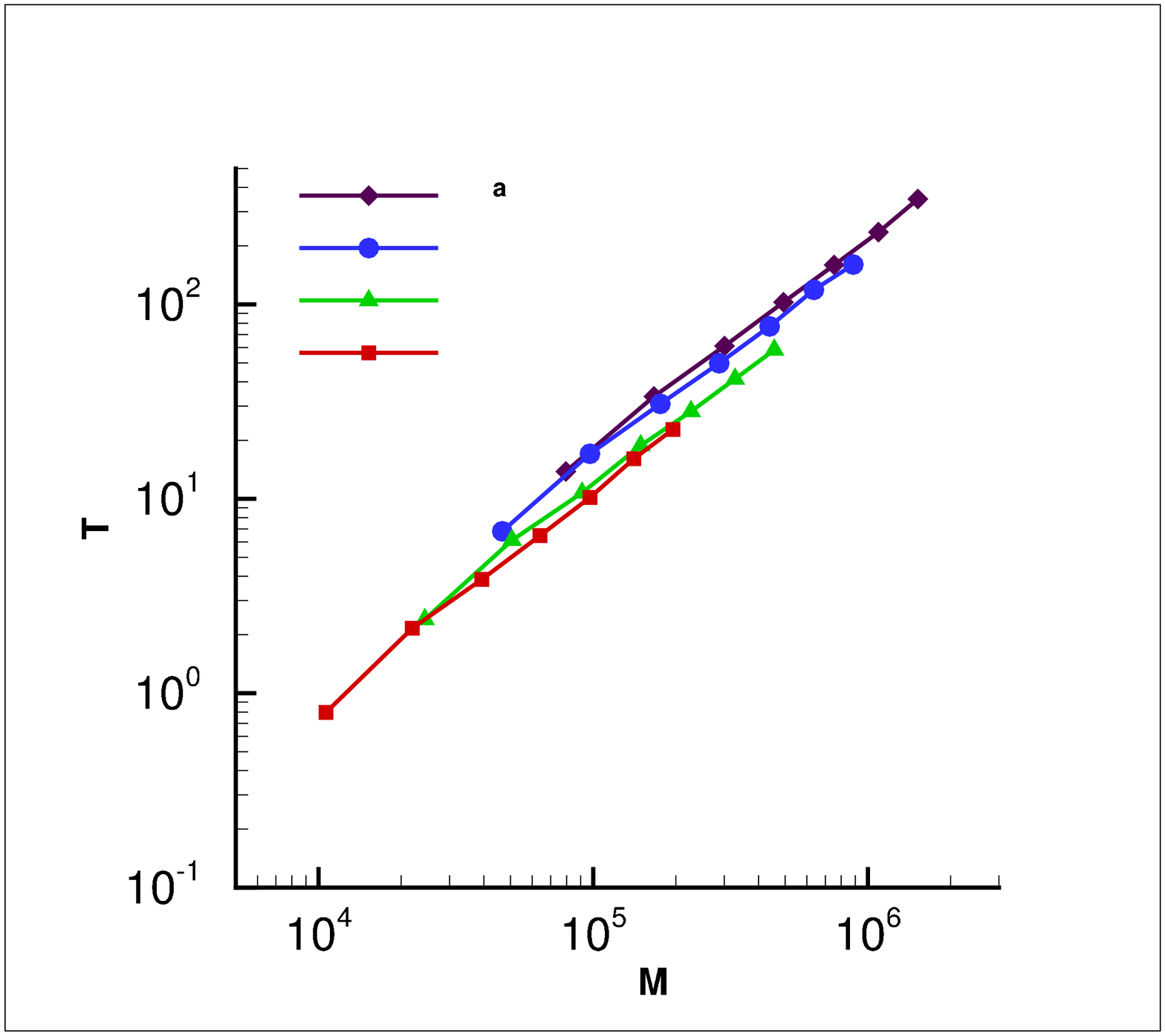}
    \caption{CPU time to solve for mesh potential versus number of mesh points $M$.  The number of elements per coordinate direction is varied from 7 to 19.}\label{fig:meshsolve}
    \end{subfigure}
    \caption{Total CPU time for (a) $N$-related components, and (b) mesh solve ($M$-related).}\label{fig:timingstotal}
\end{figure}

The computational cost of the method for $N$ charges is formulated as $\mO(N) +
\mO(M)$, where $M = {(pn_\text{el} + 1)}^3$ is the total number of degrees of
freedom in the mesh.   Example CPU time scalings for the major $N$-related
components is illustrated in Figure~\ref{fig:timingstotalN}, with the
$M$-dependent mesh solve times shown in Figure~\ref{fig:meshsolve}.  Given $N$ and
$M$ and assuming on average $\nelsr$ neighboring elements in the short-range
interaction list for each charge, then it is possible to express the
coefficients in the linear $\mO(N) + \mO(M)$ operation count in terms of $p$.
Such a formulation provides a more detailed description of the actual costs of
each component of the method and their relationships to the order of the
screens.

\subsection{Breakdown of costs}

Screens $\rho_i(\bx)$ of order $q$ are built out of ${(q+1)}^3 = \nsctotal$ basis
functions~---~recall that $p = q+2$.  The corresponding finite element solve
associated with these screens involves ${(p+1)}^3$ degrees of freedom per
element and a total of $M$ degrees of freedom.  We also define the average
number of charges per element as $\tilde{N} = N/\nel$.

\subsubsection{Screen construction}

For each evaluation, the element containing each charge is identified, and the
offsets from the center of these elements determined.  This incurs a small
$\mO(N)$ cost, which we designate $C_1 N$.  The screen coefficients are then
calculated.  As shown in (\ref{eq:moment1d}), assuming pre-computed inverses,
this amounts to three matrix-vector multiplications of size $q+1 = p-1$, for a
cost of $6{(p-1)}^2 - 3(p-1)$.  We then multiply the one-dimensional weights,
resulting in two additional floating point multiplications.  The total cost for
determining the screen coefficients is thus
\begin{equation}
 \tag{screen construction} \sim [6{(p-1)}^2 - 3(p-1) + 2{(p-1)}^3]N.
\label{eq:screencost}
\end{equation}

\subsubsection{Short-range potential}

The cost of evaluating the short-range potential depends on the method
chosen for calculating $\phsc$, as discussed in Section~\ref{sec:shortrange}.
In addition, there is a cost of $\mO(N)$ due to the singular part of the short-range calculations, which we denote $SN$.  For a general six-dimensional look-up
table, the cost of calculating $\phsc$ at a point due to all charges in the
short-range interaction volume is $C_2\tilde N\nelsr $, where $C_2$ depends on
the order of interpolation used.  If $\nsctotal$ three-dimensional look-up
tables are used, as we have done in the example calculations of Section~\ref{sec:numerics}, then the interpolation is repeated for ${(q+1)}^3$ tables and
combined by an inner product with the screen coefficients and a multiplication
by $Q_i$ for a total of $[C_2{(p-1)}^3 + 2{(p-1)}^3]\tilde N \nelsr $.  The cost
for the short-range calculation is then
\[
  \tag{point-to-point evaluation} \sim SN + [C_2{(p-1)}^3 + 2{(p-1)}^3]\tilde{N}\nelsr N.
\]

Since $\tilde{N} = N/\nel$, this expression is also
written in terms of $N^2$. However, we assume that in practice,
$\tilde{N}\nelsr$ is chosen to be small enough to render this effectively as
$\mO(N)$.  Furthermore, if $\tilde{N}$ is $> 1$, this cost is reduced
further by calculating the effects of all charges in an element at once in an
``element-to-point'' operation.  To do this, we compile a combined list of $Q_i\bc_i$
for all the charges in any given element, so that the screen potential for this sum at a
point as calculated by (\ref{eq:sc}) is the same as if the charges were
handled individually.  The cost then is then reduced by a factor of $\tilde{N}$ yielding
\[
  \tag{element-to-point evaluation} \sim SN + [C_2{(p-1)}^3 + 2{(p-1)}^3 - 1]\nelsr N.
\]

\subsubsection{Mesh solve}

The ``transfer'' of the order-$q$ screens to a representation in order $p =
q+2$ basis functions by (\ref{eq:defineassign}) to construct the source
$\rho_\text{m}$ in the right-hand side of the finite element solve
(\ref{eq:poisson}) requires an inner product between a vector containing the
screen coefficients $\bc(\bd)$ with the evaluation of the order-$q$ basis
functions at the collocation points, followed by a multiplication by $Q_i$.
This is done at each degree of freedom within an active screen area, for a
total of ${(3p+1)}^3\times[2(p-1)]N$ operations.  The multigrid solve for the
finite element problem is $\mO(M)$, with a coefficient $C$ that depends on
the convergence of the iterations, but is considered low in practice.
Overall the mesh solve thus has complexity
\[
  \tag{mesh solve} \sim \{ {(3p+1)}^3\times[2{(p-1)}^3]\}N + CM.
\]

\subsubsection{Evaluation}

The smooth potential is written as a combination of basis functions at the
location of each charge, as in (\ref{eq:evalsmooth}).  Thus evaluation
involves ${(p+1)}^3$ basis functions at a cost of $2p$ operations for each
function.  However, empirically we find that this cost is minimal in terms of
CPU time.

\subsection{Summary}

\begin{figure}[h]
  \psfrag{T}[l][1.0][1.0][-90]{$\!\!\!\!\!\!\!\!\!\!\!\!\!\!\!\!t$ (sec)}
  \psfrag{N}[c][1.0][1.0][0]{\raisebox{-0.2in}{$N$}}
  \psfrag{p}[c][1.0][1.0][0]{\raisebox{-0.1in}{$q$}}
    \psfrag{b}[c][1.0][1.0][0]{\hspace*{-0.75in}\begin{rotate}{28}\raisebox{-0.3in}{\footnotesize solid: screen potential}\end{rotate}}
      \psfrag{c}[c][1.0][1.0][0]{\hspace*{-1.2in}\begin{rotate}{28}\raisebox{-0.3in}{\footnotesize dot-dash: singular potential}\end{rotate}}
    \psfrag{d}[c][1.0][1.0][0]{\raisebox{-0.15in}{$\tilde N = 1$}}
     \psfrag{a}[][l]{\hspace*{0.73in}\raisebox{-0.62in}{\begin{minipage}{0.35\textwidth}\scriptsize$q = 4$\vspace*{0.015in}\\$q = 3$\vspace*{0.015in}\\$q = 2$\vspace*{0.015in}\\$q = 1$\end{minipage}}}
     \centering
  \begin{subfigure}[t]{0.48\textwidth}
    \includegraphics[trim=0.2in 0.05in 0.2in 0.2in, clip, width=\textwidth]{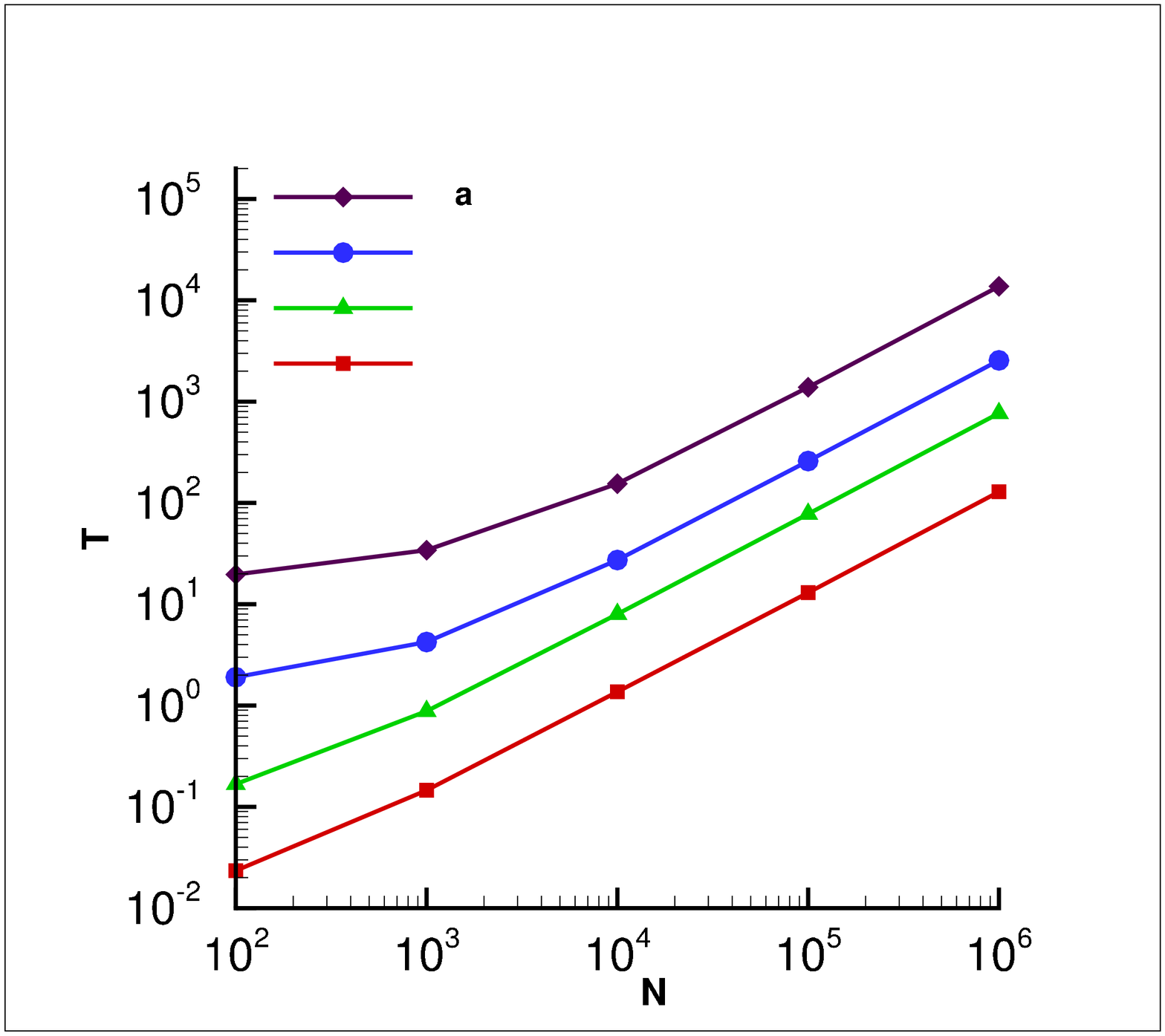}
    \caption{}
  \end{subfigure}
  \hfill
  \begin{subfigure}[t]{0.48\textwidth}
    \includegraphics[trim=0.2in 0.05in 0.2in 0.2in, clip, width=\textwidth]{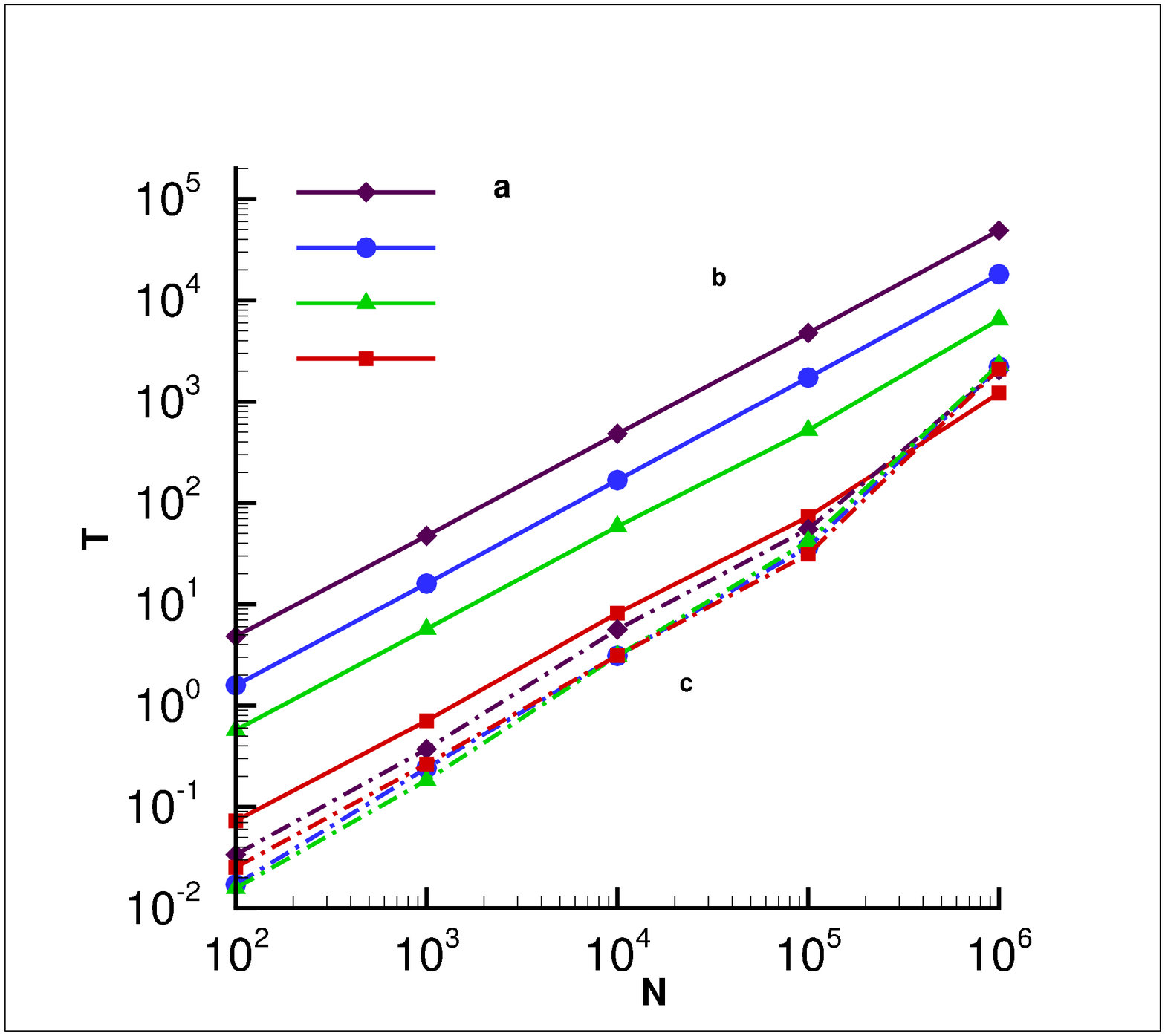}
  \caption{}
  \end{subfigure}

  \caption{Example CPU time vs.\ $N$ for the two most costly local portions of the algorithm: (a) creation of the screens, (b) calculation of short-range interactions for $N_\text{el}^\text{sr} = 7\times7\times7$.  The dot-dash lines show the time associated with calculating singular potentials, while solid lines show times calculating element-to-point screen potentials.  At large $N$, there is an expected breakdown in linearity for the singular potential calculations.}\label{fig:timingsN}
\end{figure}

The screen creation (mostly due to the ``transfer'' portion) and short-range interaction
calculations are the most costly
even for modest values of $p$ given the scaling shown above.  The relative costs
of these two portions of the algorithm depend on choices in short-range
calculation method, mesh size, and $q$.  For any given cutoff error, decreasing
mesh spacing decreases the number of short-range interactions, but results in an
increased number of collocation points $M$ in the mesh solve.  Likewise,
increasing $q$ also decreases the number of short-range interactions, but at the
price of the increased cost of constructing and manipulating screens for larger
$q$.  Calculating the short-range effects of each individual charge becomes more
costly with increased $q$, though at a slower rate than the transfer.   The
scaling of these components is shown in Figure~\ref{fig:timingsN} for cases of $N
= 10^2$ to $10^6$ randomly distributed particles in a triply-periodic box
with $6\,859$ elements.   It is noted that once $\tilde{N} \gg 1$, the singular short-range
 calculation loses its linearity in $N$.
However, the screen potential portion of the short-range calculation retains its
linearity due utilization of the ``element-to-point'' evaluation method.

\section{Example calculation}\label{sec:numerics}
\begin{figure}[h]
  \psfrag{P}[l][1.0][1.0][-90]{$\!\!\!\!\!\Phi$}
    \psfrag{Z}[l][1.0][1.0][-90]{$z$}
  \psfrag{X}[c][1.0][1.0][0]{\raisebox{-0.1in}{$x$}}
    \psfrag{x}[c][1.0][1.0][0]{\raisebox{-0.17in}{$x$}}
     \centering
  \begin{subfigure}[t]{0.48\textwidth}
    \includegraphics[trim=1.2in 0.5in 1.2in 1.in, clip, width=0.8\textwidth]{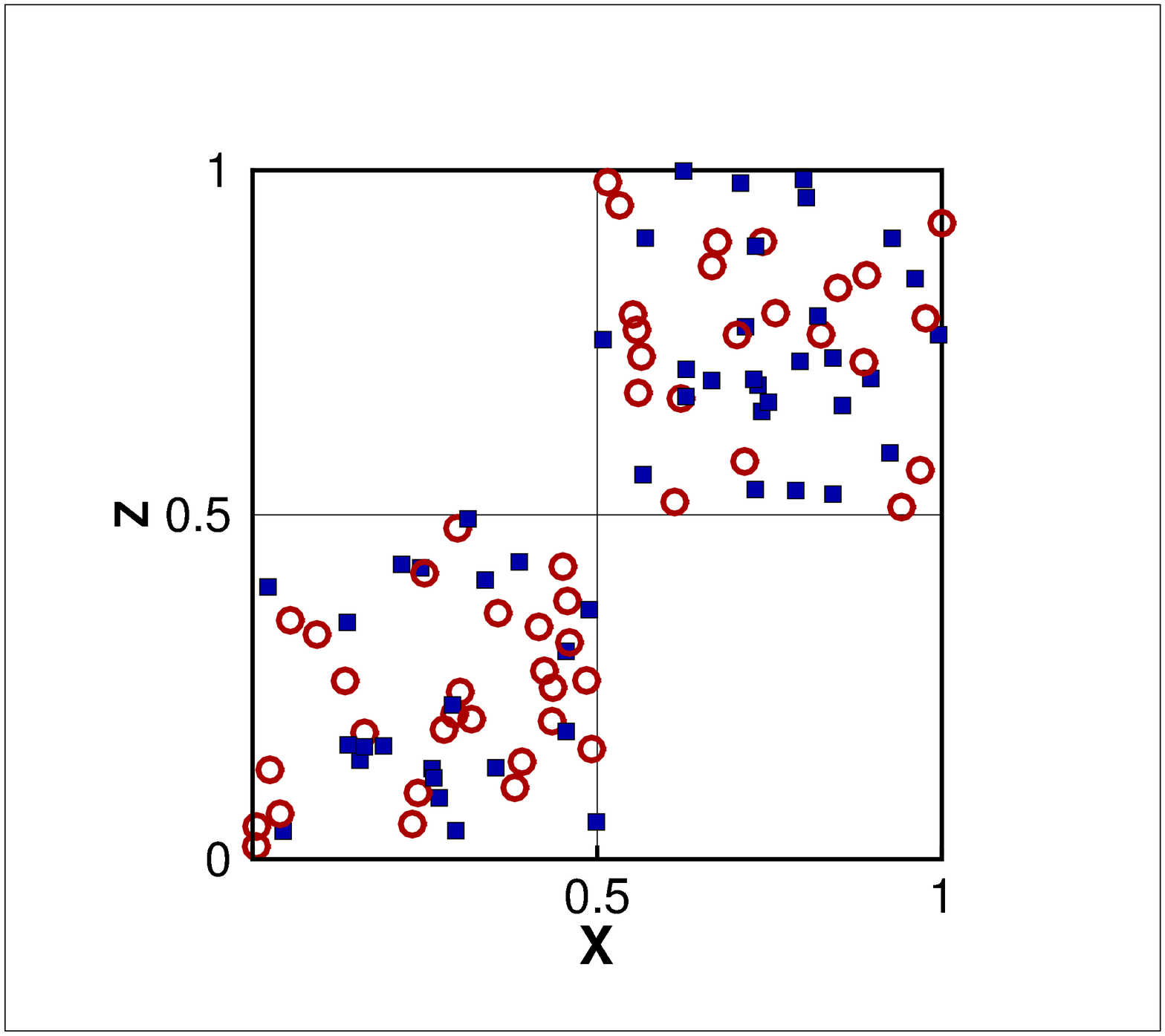}
    \caption{}
  \end{subfigure}
  \hfill
  \begin{subfigure}[t]{0.48\textwidth}
    \includegraphics[trim=0.8in 0.3in 1.in 0.5in, clip, width=0.8\textwidth]{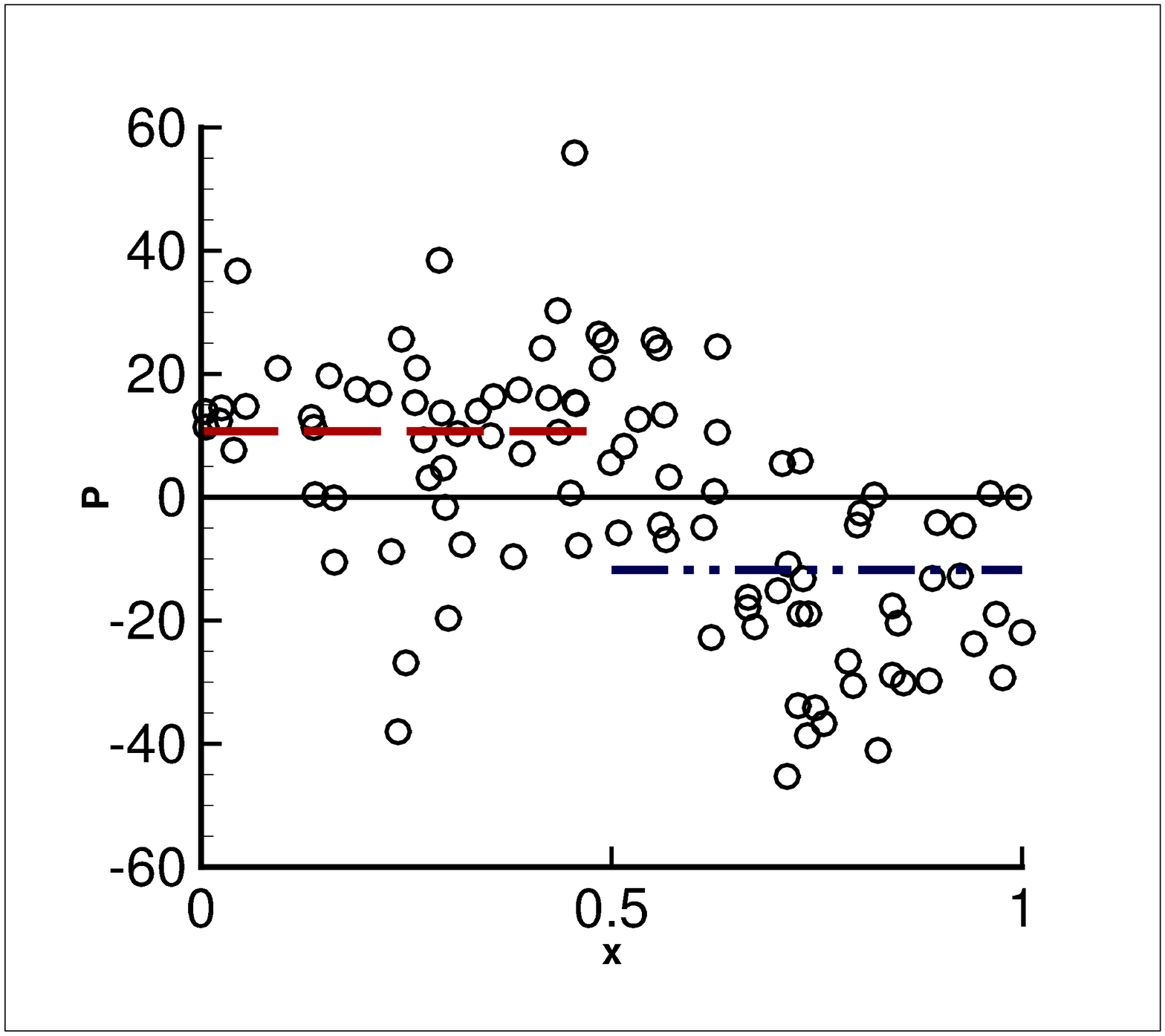}
  \caption{}
  \end{subfigure}
  \caption{Configuration of demonstration calculation: (a) Location of \color{red}$\circ$ positive \color{black} and \color{blue} $\square$ negative \color{black} charges for $N = 100$, (b) $\Phi$ at each charge location for the $100$ charge case, arranged by the charge's location in $x$; the \color{red}red \color{black}and \color{blue}blue \color{black}lines mark the average potential value at locations in the positively-biased group and negatively-biased group, respectively.}\label{fig:chargesgroups}
\end{figure}

We consider cases with  $N$ ranging from $10^2$ to $10^5$ unit charges placed
in a triply-periodic unit cube of elements with $h = 0.067$.  The exact
positions are selected randomly, but distributed so that any given charge
experiences both long-range interactions, on the scale of the overall periodic
domain size, and short-range interactions of comparable magnitude.  This is
done to provide a balanced test of both the short-range and smooth portions of
our decomposition.  To achieve this, the charges are randomly distributed
within two smaller cubes: ${[0,1/2]}^3 $ is biased toward positive charges,
$55\%$ to $45\%$, and ${[1/2,1]}^3$ is biased equally strongly toward negative
charges.  This set-up is visualized in Figure~\ref{fig:chargesgroups}a for $N =
100$.

The potential is then calculated using a short-range interaction of $7 \times 7
\times 7 = 343$ elements (corresponding to a minimum possible value of $3$ for
the cutoff distance $\hat{R}_c$) for linear through quartic screens.  This
short-range cutoff is chosen to ensure that the short-range potential of every
charge near the cutoff exhibits asymptotic behavior.   The short-range
calculation uses the approach of (\ref{eq:sc}), with $\nsctotal$ look-up
tables.  These experiments are tested using a dual, quad-core Intel Xeon E5506
CPU with 48~GB of main memory.

\begin{remark}
In our current implementation, we use a variation (but equivalent form) to this
construction, in which the values stored in the tables are for ``basis screens''
instead of screen basis functions.  These basis screens, $\rho^\text{basis}$,
are the polynomial screens associated with each node in an order-$q$ finite
element.  The values of each table are computed as a Dirichlet finite element
solution for Poisson's equation, with $-\nabla^2\Phi^\text{basis}_i =
\rho^\text{basis}_i$.  The computation is completed in a domain larger than the
size that will be kept in the look-up table to minimize boundary effects.
Because these basis screens follow our moment-canceling rules, they have
long-range decay $\sim \hat{R}^{-(q+2)}$, and the boundary conditions are accurately
set by the first terms of the multipole expansion (\ref{eq:multipole}).  The
finite element solver uses basis functions of order $p$, and the look-up tables
are stored in terms of their order-$p$ basis functions, allowing them to be
evaluated and combined in the same way as $\phsm$ for all charge locations.  We
note that because the number of tables and coefficients is unchanged, the
computational complexity for the short-range calculation is not altered by this
variation.
\end{remark}

\begin{figure}[!ht]
  \psfrag{E}[c][1.0][1.0][-90]{\hspace*{-0.8in}\text{relative error}}
  \psfrag{b}[c][1.0][1.0][0]{\raisebox{-0.5in}{\hspace*{1.0in}$\sqrt{\frac{1}{N}\sum_{i=1}^N {\left( \frac{\Phi_i-\Phi_i^\text{E}}{\Phi_i^\text{E}}\right)}^2 }$}}
    \psfrag{a}[c][1.0][1.0][0]{\raisebox{0.1in}{\hspace*{0.8in}$\max\left\{ \left| \frac{\Phi_i-\Phi_i^\text{E}}{\Phi_i^\text{E}} \right| \right\}$}}
  \psfrag{q}[c][1.0][1.0][0]{\raisebox{-0.1in}{$q$}}
  \centering
    \includegraphics[trim=0.2in 0.1in 0.2in 0.2in, clip, width=0.6\textwidth]{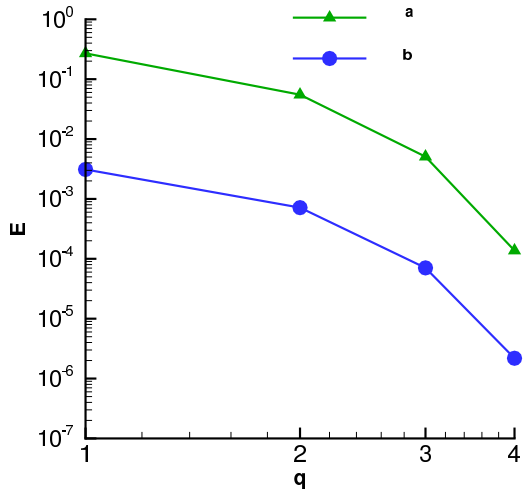}
    \caption{Convergence of the max and RMS measure of relative error in $\Phi$ versus $q$ for a sample case with $N = 10\,000$.}
\label{fig:compareerr}
\end{figure}

Upon calculation the potential is compared,
allowing for a constant which is included in a potential and in this case is equal to the average value of $\phsm$ throughout the computational domain, with that of an
Ewald summation $\Phi^\text{E}$ with large enough resolution that we consider
it the ``exact'' solution. This uses two periodic images in physical
space with $a^2 = 6.25$ and four modes for each direction in the Fourier sum.  As we see for a representative calculation in
Figure~\ref{fig:compareerr}, the method has super-algebraic convergence with
$q$.

In each of these tests, we use an algebraic multigrid preconditioned GMRES
solver with single precision residual tolerance~---~i.e., $1e-7$.
Boomeramg~\cite{boomeramg} is used and the resulting method yields 8 or fewer
iterations in each of the tests reported above.  The timing dependence on mesh size is reported in Figure~\ref{fig:meshsolve} where we see that the solver exhibits $\mathcal{O}(M)$ scaling.
   
\subsection{Estimated memory requirements}

The memory requirements of the method in our example
calculations are classified as finite element matrices or
particle-related arrays.  As the number of elements increases, the finite
element matrices comprise a majority of the total allocated memory, as demonstrated
in the following for the case of $N = 10^6$:
\begin{center}
\begin{tabular}{c | c c c}
\toprule
Order of   &  Particle       & FE matrices   & FE matrices \\
screens   &  arrays         & ($\nel = 9\times9\times9$)     & ($\nel = 25 \times 25\times25$)  \\
\midrule
linear ($q = 1$):  & 184 MB & 39 MB ($\approx$ 18\%) & 844 MB ($\approx$ 81\%) \\
quartic ($q = 4$):  & 1120 MB & 1290 MB ($\approx$ 52\%) & 27$\,$600 MB ($\approx$ 95\%) \\
\bottomrule
\end{tabular}
\end{center}
As the polynomial order increases~---~e.g.\ $q=4$, which corresponds with a
$6$th order basis for the finite element solve~---~this effect increases as
expected.  We note that this memory footprint is typical for high-order FEM, but
more optimal methods do exist~\cite{kirbybernstein2011}.

\section{Discussion}\label{sec:discussion}
\subsection{Comparison to PME}

While the method proposed here incorporates several advantageous features of
PME, there are several notable differences that may offer benefits in certain
settings.  Our method no longer relies on the FFT, which may be limiting a
extreme scales (in comparison to other Poisson solvers) and forces an assumption
of structure on the compute geometry.  The key is the introduction of a
mesh-based screen, which introduces additional complexities locally, but also
allows for a more general decomposition of the problem.  There are particle-mesh
variants that use finite elements~---~e.g., some PIC
methods~\cite{Eastwood:1986, Paes:2003}~---~but these have been proposed with a
symmetric screen, which must be resolved on the mesh.  We avoid this
approximation, but at the cost of more intricate screen functions, which are
constructed with (and the resulting potentials evaluated by) using memory-local
operations.  This fundamental difference hampers direct cost comparison with
PME/\pppm~methods, which perform well when global FFTs do not impose
restrictions.   Still, we make some general comparisons in the following.

The locality of the new method comes at the cost of more intricate screens,
which incur an $\mO(p^6)$ cost when represented by $p$-order basis functions as
discussed in Section~\ref{sec:cost}.  This is larger than the $\mO(\tilde{p}^3)$
cost of the $\tilde{p}$-order B-spline interpolations in PME.\@  However, the
polynomial order $p$ in the present scheme and the B-spline order $\tilde{p}$ in
PME are only loosely related.  The B-spline order affects the overall accuracy
of the PME method since it affects the resolution of the mesh description of
the smooth potential.  The polynomial order $p$ in the present method does not,
since the mesh solve is exact for any $p$.  Instead, $p$ affects $\hat{R}_c$
via the decay of the screened potential as shown in Figure~\ref{fig:decayavg}.
This is important, since for uniform charge density the cost of point-to-point
interactions scales with volume $\sim \hat R_c^3$.  An independent Ewald
splitting parameter sets the corresponding truncation error at fixed cut-off
radius for PME\@.

Similarly, the mesh density has different implications in the two methods.  As
with the B-spline order, the mesh density in PME affects the accuracy by
providing more resolution for the potential.  A denser mesh does not affect the
short-range calculation, but requires more global communication for the
FFT.\@  In contrast, the mesh density in the present scheme decreases
the communication burden for the short-range component of the calculation by
reducing the number of interactions included for a given $\hat{R}_c$, since the
cut-off radius is scaled by the mesh size, unlike in PME.\@  The
communication required of the mesh solver is that of multigrid.

\subsection{Comparison to FMM}

The method presented in this paper shares several
attractive features of the fast multipole method, most notably the linear
scaling.  The relative merits in comparison to FMM are likely application dependent, and
the preferred choice depends on several factors.  Though intricate,
the low communication burden of FMM leads to efficient
implementations~\cite{Lashuk:2012}.  Both methods become expensive with
increased $p$, the basis order in the present scheme or the multipole expansion
order for FMM.\@  Yet the highly local work load of the proposed high-order
screens is more suitable for emerging architectures with accelerators.  Unlike
FMM, the present method is not naturally adaptive to larger regions without
singularities~---~e.g., charges.  The degree to which FMM takes advantage
of this in parallel depends on load balancing issues of the system.

\subsection{Other considerations}

For dynamic application, the conservation properties of the overall scheme
are important, such as conservation of energy in molecular dynamics simulations.
Since we are only evaluating potentials in this paper, we do not consider
momentum or energy conservation in detail.  For the formulation as presented,
the operators we demonstrate do not exactly satisfy the symmetry discussed by
Hockney and Eastwood~\cite{Hockney:1988}, so exact momentum conservation is not
anticipated.  Moreover, as the basis functions are not differentiable at
the collocation points, straightforward analytical differentiation of the
potentials is not always possible.

The nature of the mesh solve in the presented method lends itself
to varied boundary conditions since the fundamental formulation of the
algorithm does not change when the boundary conditions are changed.  As with
PME, periodic boundary conditions are the simplest to implement in our method,
and require no extra effort beyond creating a finite element matrix that honors
the periodic structure of the mesh.  As presented in Section~\ref{sec:bc},
Dirichlet and Neumann boundary conditions simply require calculations to allow
for any short-range effects already present on the surface before applying the
conditions to the finite element problem.

The method is also extensible to non-uniform meshes common in finite element
discretizations without fundamental changes.  The main differences for general
meshes is in the cost of the method.  The screens are still built in
the same way, that is, they still solve (\ref{eq:momentmatrix}).
However, the discussed simplifications of the screen coefficient calculations
depend on a regular, rectangular mesh and are not applicable to an
unstructured mesh with general quadrilateral elements. Thus, the flexibility of
a complex mesh is balanced with the benefits of localizing the
mesh cells. Likewise, the short-range potential becomes more difficult to
generalize due to the many different shapes a screen could take based on the
shapes of the elements composing it.  Gaining accurate values for the
short-range potentials may require quadrature-based evaluations for each pair
of interacting charges.  However, the locality and structured character of
these operations is expected to coincide with high-throughput accelerators.

In our demonstration, we have presented one choice for the support of the
screens.  Another possible variation of the method is to limit the screens to
have support in only the element containing the charge, so that each screen
includes only degrees of freedom interior to the element and not those on the
faces.  Using the same multipole representation to construct the screen, this
choice results in a loss of two powers in the short-range decay of the
screens~---~e.g., $q = 3$ for the screen yields a $R^{-3}$ far-field decay
instead of $R^{-5}$, while of course still requiring $p = 5$ in the mesh solve
(and all the cost incurred by this order of $p$).  However, the more compact
screen provides an asymptotic decay rate starting at $\hat{R} \approx 1$
instead of $\hat{R} \approx 3$ (see Figure~\ref{fig:decayavg}), which
reduces the cost through reducing $\hat R_c$ for certain target accuracies.  The local
composition of these one-element screens also facilitates the move to
unstructured meshes, helping alleviate some of the additional cost in the
screen-related calculations.

In constructing our screens, we have chosen to maximize the far-field decay
rate.  Some simulation goals may be better served by other choices~---~e.g., by
a weighted objective function.  In such cases, a least-squares optimization
might provide screens with advantageous properties to meet overall simulation
objectives.  We have also restricted our discussion to purely polynomial basis
functions.  Given the regularity of the underlying Green's function, basis
enrichments with specially designed functions chosen to increase the
short-range decay likely enhance the overall performance of the method, though
this also disrupts the exactness of the mesh solve.

\section*{Acknowledgments}

This work was supported by the Computational Science \& Engineering
program at the University of Illinois at Urbana-Champaign, NSF
09--32607 and 13--36972, NSF DMS 07--46676.  This material is also based
in part upon work supported by the Department of Energy, National
Nuclear Security Administration, under Award Number DE-NA0002374.
We would also like to thank Doug Fein and the team at the National
Center for Supercomputing Applications (NCSA) for their computing
resources.

\bibliographystyle{siam.bst}
\bibliography{FErefs}

\end{document}